\documentclass[a4paper,11pt,english]{smfart}

\usepackage{amssymb}
\usepackage[T1]{fontenc}
\usepackage[cal=boondoxo]{mathalfa}

\usepackage{lscape}
\usepackage{fancybox}
\usepackage[textsize=footnotesize]{todonotes}
\usepackage[leftbars]{changebar}

\usepackage{framed}

\usepackage{palatino}
\usepackage{rotating}
\usepackage{graphicx}

\usepackage{enumerate}

\usepackage{microtype}

\usepackage{color}
\definecolor{shadecolor}{gray}{0.90}

\def\bfit{\bfseries\itshape}

\def\proofname{D\'emonstration}

\input xypic
\xyoption{all}
\xyoption{arc}

\makeindex
\def\indexnot#1#2{\index{#1@$#2$ |textbf  {\hskip0.5cm} \textbf }}
\def\indexname{Index des notations}

\newtheorem{theo}{Theorem}[section]
\def\thetheo{\thesection.\arabic{theo}}
\newtheorem{prop}[theo]{Proposition}

\newtheorem{lem}[theo]{Lemma}

\newtheorem{coro}[theo]{Corollary}

\newtheorem{defi}[theo]{Definition}

\renewcommand\theequation{\thetheo}
\def\equat{\refstepcounter{theo}\begin{equation}}
\def\endequat{\end{equation}}

\renewcommand\thetable{\thetheo}

\renewcommand\thesection{\arabic{section}}
\renewcommand\thesubsection{\thesection.\Alph{subsection}}

\def\contentsname{Table des mati\`eres}
\def\listfigurename{Liste des figures}
\def\listtablename{Liste des tables}
\def\appendixname{Appendice}

\def\refname{R\'ef\'erences}

  \def\aG{{\mathfrak a}}  
    
\def\CG{{\mathfrak C}}    \def\CM{{\mathbb{C}}}
    
  \def\eG{{\mathfrak e}}  
    
  \def\gG{{\mathfrak g}}  
    
  \def\iG{{\mathfrak i}}

\def\LG{{\mathfrak L}}  \def\lG{{\mathfrak l}}  
  \def\mG{{\mathfrak m}}  
    \def\NM{{\mathbb{N}}}
    
  \def\pG{{\mathfrak p}}  
    \def\QM{{\mathbb{Q}}}
    \def\RM{{\mathbb{R}}}
  \def\sG{{\mathfrak s}}

    \def\ZM{{\mathbb{Z}}}

  \def\ab{{\mathbf a}}  \def\AC{{\mathcal{A}}}
    
    \def\CC{{\mathcal{C}}}
    \def\DC{{\mathcal{D}}}
  \def\eb{{\mathbf e}}  \def\EC{{\mathcal{E}}}
    \def\FC{{\mathcal{F}}}
\def\Gb{{\mathbf G}}    
\def\Hb{{\mathbf H}}

  \def\kb{{\mathbf k}}  \def\KC{{\mathcal{K}}}
\def\Lb{{\mathbf L}}    \def\LC{{\mathcal{L}}}

    \def\OC{{\mathcal{O}}}
\def\Pb{{\mathbf P}}    \def\PC{{\mathcal{P}}}
    
    \def\RC{{\mathcal{R}}}
    
  \def\tb{{\mathbf t}}  
  \def\ub{{\mathbf u}}

\def\Xb{{\mathbf X}}    \def\XC{{\mathcal{X}}}
    
\def\Zb{{\mathbf Z}}  \def\zb{{\mathbf z}}  \def\ZC{{\mathcal{Z}}}

\def\Crm{{\mathrm{C}}}    
    
\def\Erm{{\mathrm{E}}}    
\def\Frm{{\mathrm{F}}}

\def\Mrm{{\mathrm{M}}}    
\def\Nrm{{\mathrm{N}}}    
    
  \def\prm{{\mathrm{p}}}

\def\Srm{{\mathrm{S}}}  \def\srm{{\mathrm{s}}}  
\def\Trm{{\mathrm{T}}}

\def\Zrm{{\mathrm{Z}}}

\def\a{\alpha}

\def\g{\gamma}
\def\G{\Gamma}
\def\d{\delta}
\def\D{\Delta}
\def\e{\varepsilon}
\def\ph{\varphi}
\def\l{\lambda}
\def\L{\Lambda}

\def\o{\omega}
\def\O{\Omega}
\def\r{\rho}

\def\t{\tau}

\def\z{\zeta}

\def\lamb{{\boldsymbol{\lambda}}}       
       
\def\mub{{\boldsymbol{\mu}}}

\def\Omeb{{\boldsymbol{\Omega}}}        
            \def\pit{{\tilde{\pi}}}
            \def\Pit{{\tilde{\Pi}}}

\def\pibt{{\tilde{\boldsymbol{\pi}}}}

\DeclareMathOperator{\End}{{\mathrm{End}}}

\DeclareMathOperator{\Hom}{{\mathrm{Hom}}}

\DeclareMathOperator{\Id}{{\mathrm{Id}}}

\DeclareMathOperator{\Irr}{{\mathrm{Irr}}}
\DeclareMathOperator{\Ker}{{\mathrm{Ker}}}

\DeclareMathOperator{\rad}{{\mathrm{rad}}}
\DeclareMathOperator{\Rad}{{\mathrm{Rad}}}
\DeclareMathOperator{\Res}{{\mathrm{Res}}}

\DeclareMathOperator{\Pic}{{\mathrm{Pic}}}
\DeclareMathOperator{\Movable}{{\mathrm{Mov}}}

\def\to{\rightarrow}
\def\longto{\longrightarrow}
\def\injto{\hookrightarrow}

\def\fonction#1#2#3#4#5{\begin{array}{rccc}
{#1} : & {#2} & \longto & {#3}  \\
& {#4} & \longmapsto & {#5}
\end{array}}

\def\fonctio#1#2#3#4{\begin{array}{ccc}
{#1} & \longto & {#2} \\
{#3} & \longmapsto & {#4}
\end{array}}

\def\DS{\displaystyle}

\def\SSS{\scriptscriptstyle}

\def\finl{~$\blacksquare$}

\def\lexp#1#2{\kern\scriptspace\vphantom{#2}^{#1}\kern-\scriptspace#2}
\def\le{\hspace{0.1em}\mathop{\leqslant}\nolimits\hspace{0.1em}}
\def\ge{\hspace{0.1em}\mathop{\geqslant}\nolimits\hspace{0.1em}}

\mathchardef\inferieur="321E
\mathchardef\superieur="321F

\def\eqna{\begin{eqnarray*}}
\def\endeqna{\end{eqnarray*}}

\def\mini{{\mathrm{min}}}
\def\itemth#1{\item[${\mathrm{(#1)}}$]}

\catcode`\@=11
\long\def\@car#1#2\@nil{#1}
\long\def\@first#1#2{#1}
\long\def\@second#1#2{#2}
\long\def\ifempty#1{\expandafter\ifx\@car#1@\@nil @\@empty
  \expandafter\@first\else\expandafter\@second\fi}
\catcode`\@=12

\def\GL{{\mathrm{GL}}}

\DeclareMathOperator{\REF}{Ref}

\def\boitegrise#1#2{\begin{centerline}{\fcolorbox{black}{shadecolor}{~
    \begin{minipage}[t]{#2}{\vphantom{~}#1\vphantom{$A_{\DS{A_A}}$}}
            \end{minipage}~}}\end{centerline}\medskip}

\def\ve{{\SSS{\vee}}}

\def\surto{\twoheadrightarrow}

\theoremstyle{remark}
\newtheorem{rema}[theo]{Remark}

\newtheorem{exemple}[theo]{Example}

\theoremstyle{plain}

\newtheorem{conjecture}[theo]{Conjecture}

\def\Frac{{\mathrm{Frac}}}

\def\BIL{LR}
\def\GAUCHE{L}
\def\CAR{CAR}
\def\FAM{FAM}

\def\groth{\KC_0}

\def\reg{{\mathrm{reg}}}

\def\refw{{\REF(W)/W}}

\def\euler{{\eb\ub}}

\def\calo{{\Crm\Mrm}}

\def\xyinj{\ar@{^{(}->}}
\def\xysur{\ar@{->>}}

\def\isomorphisme#1{{\boldsymbol{[}}\hskip0.5mm #1\hskip0.5mm {\boldsymbol{]}}}

\DeclareMathOperator{\carac}{{\mathrm{Car}}}
\bigskip

\makeatletter
\def\hlinewd#1{%
\noalign{\ifnum0=`}\fi\hrule \@height #1 %
\futurelet\reserved@a\@xhline}
\makeatother

\newlength\epaisLigne

\usepackage{multirow}

\newcommand{\longiso}{\stackrel{\sim}{\longrightarrow}}

\def\LCov{{\overline{\LC}}}

\def\carac{{\mathrm{car}}}

\makeatletter
\def\hlinewd#1{%
\noalign{\ifnum0=`}\fi\hrule \@height #1 %
\futurelet\reserved@a\@xhline}
\makeatother

\usepackage{multirow}

\usepackage{slashbox}

\usepackage{lscape}

\def\troncation{{\mathrm{Trunc}}}

\usepackage{pstricks}

\usepackage{listings}  
\begin{document}


\title{Computational aspects of Calogero--Moser spaces}

\author{{\sc C\'edric Bonnaf\'e}}
\address{IMAG, Universit\'e de Montpellier, CNRS, Montpellier, France} 
\email{cedric.bonnafe@umontpellier.fr}

\author{{\sc Ulrich Thiel}}

\address{
Department of Mathematics, University of Kaiserslautern, Postfach 3049, 67653 Kaiserslautern, Germany}
\email{thiel@mathematik.uni-kl.de}

\date{\today}

\thanks{The first author is partly supported by the ANR:
Projects No ANR-16-CE40-0010-01 (GeRepMod) and ANR-18-CE40-0024-02 (CATORE).
This work is a contribution to the SFB-TRR 195 'Symbolic Tools in Mathematics and their Application' of the German Research Foundation (DFG). We would like to thank Gwyn 
Bellamy and Gunter Malle for comments on a preliminary version of this paper.}

\maketitle
\pagestyle{myheadings}

\markboth{\sc C. Bonnaf\'e \& U. Thiel}{\sc Computational aspects of Calogero--Moser spaces}

\begin{abstract}
We present a series of algorithms for computing geometric and representation-theoretic invariants of Calogero--Moser spaces and rational Cherednik algebras associated with complex reflection groups. In particular, we are concerned with Calogero--Moser families (which correspond to the $\mathbb{C}^\times$-fixed points of the Calogero--Moser space) and cellular characters (a proposed generalization by Rouquier and the first author of Lusztig's constructible characters based on a Galois covering of the Calogero--Moser space). To compute the former, we devised an algorithm for determining generators of the center of the rational Cherednik algebra (this algorithm has several further applications), and to compute the latter we developed an algorithmic approach to the construction of cellular characters via Gaudin operators. We have implemented all our algorithms in the Cherednik Algebra Magma Package (CHAMP) by the second author and used this to confirm open conjectures in several new cases. 
As an interesting application in birational geometry we are able to determine for many exceptional complex reflection groups the chamber decomposition of the movable cone of a $\mathbb{Q}$-factorial terminalization (and thus the number of non-isomorphic relative minimal models) of the associated symplectic singularity. 
Making possible these computations was also a source of inspiration for the first author to propose 
several conjectures about the geometry of Calogero--Moser spaces (cohomology, fixed points, 
symplectic leaves), often in relation with the representation theory of finite reductive 
groups.
\end{abstract}

\section{Introduction}
Let $V$ be a finite-dimensional complex vector space and let $W \subseteq \GL(V)$ be a finite subgroup generated by reflections in $V$. The {\it rational Cherednik algebras} defined by Etingof and Ginzburg \cite{EG} yield a flat family of deformations $\Hb_c$ of the algebra $\mathbb{C} \lbrack V \times V^* \rbrack \rtimes W$. Here, $\mathbb{C} \lbrack V \times V^* \rbrack$ is the coordinate ring of the space $V \times V^*$ considered as an algebraic variety and $\mathbb{C} \lbrack V \times V^* \rbrack \rtimes W$ denotes the semi-direct product with the group algebra of~$W$, which incorporates the action of $W$ on $V \times V^*$. Moreover, $c$ is a parameter from a complex vector space $\CC$ whose dimension is equal to the number of conjugacy classes of reflections contained in $W$. We note that in \cite{EG} a further deformation parameter $t \in \mathbb{C}$ is considered but here we focus on the case $t=0$. Since their introduction, rational Cherednik algebras have attracted a tremendous amount 
of interest. This stems from their diverse connections to other fields and problems, especially in algebraic geometry and representation theory. We will mention two such connections that serve as motivation for this paper.

The spectrum $\ZC_c$ of the center $\Zb_c$ of $\Hb_c$ is the {\it Calogero--Moser space} originating from physics \cite{CalogeroOriginal, MoserOriginal}. It is a Poisson deformation of the quotient variety $(V \times V^*)/W$, the latter being a symplectic singularity in the sense of Beauville \cite{beauville}. The representation theory of $\Hb_c$ is closely tied to the geometry of $\ZC_c$ and plays a key role in the question of whether $(V \times V^*)/W$ admits a symplectic (equivalently, crepant) resolution, see \cite{EG,GK,Namikawa-Poisson,BST-Towards}, and in determining the chamber decomposition of the movable cone of a $\mathbb{Q}$-factorial terminalization (i.e.\ a relative minimal model) of $(V \times V^*)/W$, see \cite{Namikawa-Birational,BST-Hyperplanes}.

The representation theory of $\Hb_c$ is furthermore (conjecturally) linked to (parts of) Kazhdan--Lusztig theory \cite{KL} and provides a candidate for extending this theory from finite Coxeter groups to complex reflection groups following the philosophy of the ``spetses'' program \cite{BMM}. It was noticed in \cite{gordon} that $\Zb_c$ contains the subalgebra $\Pb = \CM[V]^W \otimes \CM[V^*]^W$. This inclusion defines a finite $\CM^\times$-equivariant morphism $\Upsilon_c \colon \ZC_c \to \PC$. The quotient of $\Hb_c$ by the ideal generated by the origin in $\PC$ is the {\it restricted rational Cherednik algebra} $\overline{\Hb}_c$. Its blocks are in bijection with the fiber of $\Upsilon_c$ at the origin---which in turn is precisely the set of $\CM^\times$-fixed points of $\ZC_c$---and the set of simple modules of $\overline{\Hb}_c$ is naturally in bijection with the set $\Irr(W)$ of irreducible complex characters of $W$. In particular, the blocks of $\overline{\Hb}_c$ yield a partition of $\Irr(W)$. These are the {\it Calogero--Moser families}. They are (conjecturally) related to Rouquier families defined by the associated Hecke algebra, see \cite{gordon martino, martino, chlouveraki LNM}, which generalize Lusztig's families of unipotent characters associated with finite Coxeter groups \cite{lusztig orange,lusztig}. Rouquier and the first author \cite{calogero-first, calogero} noticed that the Galois group of
the Galois
closure $\RC_c \to \PC$ of the covering $\ZC_c \to \PC$ acts
(non-canonically) on the set $W$. By considering orbits in $W$ under
various subgroups of this Galois group, they obtained a decomposition of
$W$. These are the {\it Calogero--Moser cells} which are conjectured to coincide with Kazhdan--Lusztig cells and thus potentially provides a way to generalize them to complex reflection groups. This conjecture was recently proven to be true for the symmetric group \cite{cells-typeA}. Associated to each Calogero--Moser left cell \cite{calogero} constructed a {\it cellular character}: several descriptions
are available and, in this paper, we will use the one involving Gaudin operators.
They are expected to coincide with Lusztig's constructible characters~\cite{lusztig}
for Coxeter groups.

\bigskip
\subsection{Algorithms}

According to the Shephard--Todd classification \cite{ST}, the class of irreducible complex reflection groups splits into two classes: the infinite series of groups $G(d,m,n)$, which are normal subgroups of the wreath product $C_d \wr S_n$, and 34 exceptional groups $G_4, \ldots, G_{37}$. Whereas the  infinite series shows a common combinatorial pattern, each of the exceptional groups is essentially unique and often needs separate treatment.

It is especially because of the exceptional groups that we wanted to find algorithms allowing us to explicitly compute invariants of the associated Calogero--Moser spaces and rational Cherednik algebras so that we can test conjectures and gather data to drive the development of the theory. An important aspect of our approach is that we wanted to compute invariants for {\it all} parameters $c$. This is why we usually work with the {\it generic} rational Cherednik algebra $\Hb$ (and its center $\Zb$) defined over a polynomial ring, i.e.\ our parameters $c$ are indeterminates.

In this paper we present algorithms for computing the following:
\begin{enumerate}
\item The inverse of the natural truncation map $\troncation \colon \Zb \to \CM[\CC \times V \times V^*]^W$, which is an isomorphism of $\CM[\CC]$-modules (Algorithm \ref{algo:troncation}). Our algorithm iteratively deforms an element of $\CM[\CC \times V \times V^*]^W$ to an element of $\Zb$.
\item A minimal system of algebra generators of $\Zb$ (Algorithm \ref{algo:generateurs}).
\item A presentation of $\Zb$ and thus of the Calogero--Moser space (Algorithm \ref{algo:relations}).
\item The Calogero--Moser families (Algorithm \ref{algo:familles}).
\item The Calogero--Moser hyperplanes (Algorithm \ref{algo:familles}). The complement of this hyperplane arrangement is precisely the locus where the Calogero--Moser families remain {\it generic}. These hyperplanes play a similar role as the essential hyperplanes for Rouquier families \cite{chlouveraki LNM}.
\item Cuspidal families (Algorithm \ref{algo:cuspidal}). These are the families lying on a zero-dimensional symplectic leaf of $\ZC_c$, see \cite{bellamy cuspidal, bellamy thiel}.
\item Cellular characters (Algorithm \ref{algo:cellulaires}). Our approach uses the description of cellular characters via Gaudin operators from \cite{calogero}.
\end{enumerate}

\bigskip
\subsection{Results and applications}
Our algorithms are not just theoretical. We have implemented all of them in the Cherednik Algebra Magma Package (CHAMP) created by the second author \cite{champ}, which is based on the computer algebra system MAGMA \cite{magma} and is available at
\begin{center}
 {\tt https://github.com/ulthiel/Champ} \;.
\end{center}
Because of the complexity and wealth of data—there are many invariants and all of them depend on parameters—we can in this paper only give a brief summary of what is now actually computable. Some of the data is organized in a database which is accessible from within CHAMP. The most striking results and applications presented 
in this paper are the following:
\begin{enumerate}
\item An explicit minimal system of algebra generators of $\Zb$ is now known for all exceptional complex reflection groups except $G_{16} - G_{22}$ and $G_{27} - G_{37}$; moreover, all algebra generators of $\mathbb{Z}$-degree 0 are known also for $G_{27}$ and $G_{28}$ (see Section \ref{center_comp_overview}). Knowing the degree-0 generators is the key to the following results on families and hyperplanes.

\item The Calogero--Moser hyperplanes are now known for all exceptional complex reflection groups except $G_{16} - G_{19}$, $G_{21}$, and $G_{32}$ (see Table \ref{table:cm_hyperplanes}). Applying results of~\cite{Namikawa-Birational,BST-Hyperplanes} we thus know in all these cases the chamber decomposition of the movable cone of a $\mathbb{Q}$-factorial terminalization (i.e.\ a relative minimal model) of the associated symplectic singularity $(V \times V^*)/W$ and we know the number of non-isomorphic relative minimal models. Particularly exciting is that we were able to determine this for the Weyl group $G_{28} = F_4$ (the computation took several days).

\item The Calogero--Moser families are now known for all exceptional complex reflection groups (and {\it all} parameters) except $G_{16} - G_{19}$, $G_{21}$, $G_{29}-G_{37}$. In all cases we confirmed a conjecture by Martino \cite{martino} stating that the Calogero--Moser families are unions of the Rouquier families
(Theorem~\ref{theo:gordon-martino}). Particularly exciting is again the case $G_{28}=F_4$ where the Rouquier families are the same as the Lusztig families and the equality to Calogero--Moser families was conjectured previously by Gordon and Martino \cite{gordon martino}. In this last (big) case, we also prove
that there is a unique cuspidal Calogero--Moser family, confirming a conjecture by Bellamy
and the second author~\cite{bellamy thiel}.

\item The Calogero--Moser cellular characters are now known for {\it spetsial} groups of rank $\le 2$
at {\it spetsial} parameters (see~\cite{malle icm} for the definition of spetsial): in all known cases, they coincide with the {\it cellular characters} defined by Malle and
Rouquier~\cite{malle rouquier} built from the {\it spetses} philosophy~\cite{BMM}.
\end{enumerate}
Apart from the above results, making possible these 
computations was used in several other papers where the first author is involved, 
to understand some symplectic singularities~\cite{BBFJLS} or to propose 
several conjectures about the geometry of Calogero--Moser spaces (cohomology, fixed points, 
symplectic leaves), often in relation with the representation theory of finite reductive 
groups~\cite{auto, regular, cm-unip}. Let us mention some of these applications:
\begin{enumerate}
\item For dihedral groups of order $8$ and $12$ the symplectic singularity in the origin of the associated Calogero--Moser space has been identified in~\cite{BBFJLS} as the symplectic singularity of the closure of the minimal orbit of a semisimple Lie algebra (see Proposition \ref{prop:b2g2}). Basis for this identification was the computation of $\Zb$ via our algorithm, together with its Poisson bracket. We were now able to achieve the identification also in case of the exceptional group $G_4$ (Proposition \ref{eq:g4-sl3}).

\item Explicit equations of the Calogero--Moser spaces associated with 
dihedral groups at equal parameters have been obtained in~\cite{bonnafe diedral 2}, 
together with the explicit Poisson bracket. The pattern was found thanks to 
experiments for dihedral groups of order $8$, $10$, $12$, $14$. 
This was applied by Bellamy, Fu, Juteau, Levy, Sommers and 
the first author to find a new family of isolated symplectic singularities with trivial 
local fundamental group~\cite{BBFJLS}, answering an old question of Beauville~\cite{beauville}. 

\item It is conjectured in~\cite[Conj.~B]{auto} that the normalization of the closure of a symplectic 
leaf in the fixed point subvariety of the Calogero--Moser space under some automorphism is also 
a Calogero--Moser space (associated with another complex reflection group, perfectly identified). 
For $G_4$ we could confirm this conjecture (Theorem \ref{theo:g4}).

\item The first author proposed an intriguing conjecture~\cite[Conj.~12.3]{cm-unip} 
relating the parameters involved in the previous conjecture~\cite[Conj.~B]{auto} 
with the parameters involved in the endomorphism algebra of the cohomology of a 
Deligne-Lusztig variety. This was inspired by experiments in type $B_2$ or $G_2$ 
and, for people believing in {\it spetses} (see~\S\ref{sub:spetses}), 
by explicit computations for $G_4$ (see~\cite[\S{18}]{cm-unip}). 
\end{enumerate}

\tableofcontents

\def\carac{{\mathrm{char}}}
\def\rad{{\mathrm{rad}}}
\def\semi{{\mathrm{sem}}}

\def\red{{\mathrm{red}}}
\def\nor{{\mathrm{nor}}}
\def\setw{{\mathrm{set}}}
\def\ptw{{\mathrm{pt}}}
\def\la{\langle}
\def\ra{\rangle}

\bigskip
\section{Calogero--Moser spaces and rational Cherednik algebras} \label{sec:setup}

We recall some basics about Calogero--Moser spaces and rational Cherednik algebras from \cite{EG}, keeping the notation of~\cite{calogero}. Throughout this paper, we abbreviate $\otimes_\CM$ by~$\otimes$. All varieties are complex algebraic, and for an affine variety $\XC$ we denote its coordinate ring by $\CM[\XC]$. For a group $W$ we denote by $\Irr(W)$ the set of complex irreducible characters of $W$.


\medskip
\def\alephb{{\boldsymbol{\aleph}}}

\subsection{Reflection groups}
We fix in this paper a complex vector space $V$ of finite dimension $n$ and
%
%
a finite
subgroup $W$ of $\GL(V)$. For an element $w \in W$ we denote by $V^w$ the fixed space of $w$. By
$$\REF(W)=\{s \in W~|~\dim_\CM V^s=n-1\}$$
we denote the set of {\it reflections} in $W$. We assume throughout this paper that $W$ is generated by its reflections, i.e. $W$ is a {\it reflection group}.
We denote by $$\e : W \to \CM^\times \;, \; w \mapsto \det(w) \;,$$ the determinant character of $W$.
Considering $V$ and $V^*$ as algebraic varieties, we identify the coordinate ring $\CM[V]$ (resp. $\CM[V^*]$) of $V$ (resp. $V^*)$ with the symmetric
algebra $\Srm(V^*)$ (resp. $\Srm(V)$).
If $s \in \REF(W)$, we denote by $\a_s \in V^*$ (resp. $\a_s^\ve \in V$) an element
such that $(V^s)^\perp=\CM \a_s$ (resp. $(V^{*s})^\perp=\CM \a_s^\ve$). Here, $(V^s)^\perp$ denotes the space of linear forms whose kernel contains $V^s$, and $(V^{*s})^\perp$ is defined analogously.

We denote by $\AC$ the set of {\it reflecting hyperplanes} of $W$, namely
$$\AC=\{V^s~|~s \in \REF(W)\} \;.$$
If $H \in \AC$, we denote by $W_H$ the pointwise stabilizer of $H$, by
$\a_H$ an element of $V^*$ such that
$H=\Ker(\a_H)$, and by $\a_H^\vee$ an element such that
$V=H \oplus \CM \a_H^\vee$ and the line $\CM\a_H^\vee$ is $W_H$-stable. For $s \in \REF(W)$,
we set $\a_s=\a_{V^s}$ and $\a_s^\vee=\a_{V^s}^\vee$. 
We set $e_H=|W_H|$. Note that $W_H$ is cyclic of order $e_H$ and that
$$\Irr(W_H)=\{\Res_{W_H}^W \e^j~|~0 \le j \le e_H-1\} \;.$$ We denote by $\e_{H,j}$
the (central) primitive idempotent of $\CM W_H$ associated with the character
$\Res_{W_H}^W \e^{-j}$, namely
$$\e_{H,j}=\frac{1}{e_H}\sum_{w \in W_H} \e(w)^j w \in \CM W_H.$$
If $\O$ is a $W$-orbit of reflecting hyperplanes, we write $e_\O$ for the
common value of all the $e_H$, where $H \in \O$.
We denote by $\aleph(W)$ the set of pairs $(\O,j)$ where $\O \in \AC/W$ and
$0 \le j \le e_\O-1$.

\bigskip

\subsection{Parameters} \label{sec:parameters}
Let $\CC$ denote the $\CM$-vector space of maps $$c : \REF(W) \longto \CM \;, \;
s \mapsto c_s \;,$$ which are invariant under $W$-conjugation. We denote by $\REF(W)/\!\sim$ the
set of conjugacy classes of reflections of $W$. For $s \in \REF(W)$
we denote by $C_s \in \CC^*$ the linear form defined by $C_s(c)=c_s$. Of course, $C_s=C_t$
if $s$ and $t$ are conjugate. We identify $\CM[\CC]=\CM[(C_s)_{s \in \REF(W)/\!\sim}]$.
If $c \in \CC$, we denote by $\CG_c$ the ideal of $\CM[\CC]$ of functions
vanishing at $c$. This ideal is generated by $(C_s-c_s)_{s \in \REF(W)/\!\sim}$.

Let $\KC$ denote the space of maps $$k : \aleph(W) \longto \CM \;, \; (\O,j) \mapsto k_{\O,j} \;,$$
which satisfy $k_{\O,0}+k_{\O,1}+\cdots + k_{\O,e_\O-1}=0$ for all $\O \in \AC/W$.
For $k \in \KC$, $H \in \AC$, and $0 \le j \le e_H-1$ we write
$k_{H,j}=k_{\O,j}$,
where $\O$ is the $W$-orbit of $H$. The linear
map
$$\g : \KC \to \CC \;, \; \g(k)_s=\sum_{j=0}^{e_{V^s}-1} \e(s)^{j-1} k_{V^s,j} \;,$$
is an isomorphism and we denote by $$\kappa : \CC \to \KC$$ its inverse.
Throughout the paper, we will identify $\CC$ and $\KC$ through these
isomorphisms.

Let $K_{\O,j} \in \KC^*$ be the linear form defined by $K_{\O,j}(k)=k_{\O,j}$ for all $k \in \KC$.
Then $K_{\O,0}+K_{\O,1}+\cdots + K_{\O,e_\O-1}=0$. As before, we write $K_{H,j}=K_{\O,j}$. Through the
identification $\KC^* \simeq \CC^*$, we have
\equat\label{eq:k-c}
C_s=\sum_{j=0}^{e_{V^s}-1} \e(s)^{j-1} K_{V^s,j}.
\endequat
Note that $(K_{(\O,j)})_{(\O,j) \in \aleph(W), j \neq 0}$ is a basis of $\KC^*=\CC^*$.

\bigskip

\subsection{Rational Cherednik algebras at ${\boldsymbol{t=0}}$}
We follow the convention of~\cite{calogero} and
define the {\it generic rational Cherednik algebra $\Hb$} ({\it at $t=0$}) to be the quotient
of the $\CM[\CC]$-algebra $\CM[\CC] \otimes (\Trm(V\times V^*)\rtimes W)$, the second factor being the semi-direct product of the tensor algebra $\Trm(V \times V^*)$ with the group $W$, by the relations
\equat\label{eq:rels}
\begin{cases}
[x,x']=[y,y']=0,\\
[y,x]=\DS{\sum_{s \in \REF(W)} (\e(s)-1)C_s
\frac{\langle y,\a_s \rangle \cdot \langle \a_s^\ve,x\rangle}{\la \a_s^\ve,\a_s\ra}\, s},
\end{cases}
\endequat
for all $x$, $x'\in V^*$, $y$, $y'\in V$. Here $\la\ ,\ \ra: V\times V^*\to\CM$ is the standard pairing.
From~\eqref{eq:k-c}, one deduces an
equivalent presentation, namely
\equat\label{eq:rels-k}
\begin{cases}
[x,x']=[y,y']=0,\\
[y,x]=\DS{\sum_{H\in\mathcal{A}} \sum_{j=0}^{e_H-1}
e_H(H_{H,j}-K_{H,j+1})
\frac{\langle y,\a_H \rangle \cdot \langle \a_H^\ve,x\rangle}{\langle \a_H^\ve,\a_H\rangle} \e_{H,j}},
\end{cases}
\endequat
for all $x$, $x'\in V^*$, $y$, $y'\in V$.

The first commutation relations imply that
we have morphisms of algebras $\CM[V] \to \Hb$ and $\CM[V^*] \to \Hb$.
By~\cite[Theo.~1.3]{EG}
 we have an isomorphism of $\CM$-vector spaces
\equat\label{eq:pbw}
\CM[\CC] \otimes \CM[V] \otimes \CM W \otimes \CM[V^*] \longiso \Hb
\endequat
induced by multiplication. This is the so-called {\it PBW-decomposition} of $\Hb$.

%
%

\bigskip

\subsection{Calogero--Moser space}
We denote by $\Zb$ the center of the algebra $\Hb$. It is well-known~\cite{EG} that
$\Zb$ is an integral domain which is integrally closed. Moreover, $\Zb$ contains
$\CM[V]^W$ and $\CM[V^*]^W$ as subalgebras~\cite{gordon}, hence it contains $$\Pb=\CM[\CC] \otimes \CM[V]^W \otimes \CM[V^*]^W \;,$$
and $\Zb$ is a free $\Pb$-module of rank $|W|$. We denote by $\ZC$ the
affine algebraic variety whose ring of regular functions $\CM[\ZC]$ is $\Zb$.
This is the {\it generic Calogero--Moser space} associated with $(V,W)$.
It is irreducible and normal.

We set $$\PC=\CC \times V/W \times V^*/W \;,$$ so $\CM[\PC]=\Pb$. The inclusion
$\Pb \injto \Zb$ induces a morphism of varieties
\equat \label{upsilon}
\Upsilon : \ZC \longto \PC
\endequat
which is finite and flat.

For $c \in \CC$ we denote by $$\Hb_c=\Hb/\CG_c \Hb$$ the {\it specialization} of $\Hb$ at $c$. Similarly, we write $\Zb_c=\Zb/\CG_c \Zb$ and $\Pb_c = \Pb/\CG_c \Pb$. It turns out~\cite[Coro.~4.2.7]{calogero} that $\Zb_c$ is the
center of $\Hb_c$. We let $\ZC_c$ denote the variety $\Upsilon^{-1}(\{c\} \times V/W \times V^*/W)$,
so $\CM[\ZC_c]=\Zb_c$.


\bigskip

%

\subsection{Extra-structures on the center}
The center $\Zb$ is endowed with extra-structures
(a bi-grading, a Poisson bracket, and a special central element)
which are described below.

\bigskip

\subsubsection{Gradings and $\CM^\times$-action}
The $\CM$-algebra $\CM[\CC] \otimes (\Trm(V\times V^*)\rtimes W)$ can be $(\NM \times \NM)$-graded
in such a way that the generators have the following bi-degrees
$$
\begin{cases}
\deg(C_s) =(1,1) & \text{if $s \in \REF(W)/\!\sim$,}\\
\deg(y)=(0,1) & \text{if $y \in V$,}\\
\deg(x)=(1,0) & \text{if $x \in V^*$,}\\
\deg(w)=0 & \text{if $w \in W$.}
\end{cases}
$$
This descends to an $(\NM \times \NM)$-grading on $\Hb$
because the defining relations~(\ref{eq:rels})
are bi-homogeneous. As a consequence, $\Zb$ is also $(\NM \times \NM)$-graded.

This bi-grading induces a $\ZM$-grading on $\Hb$ and $\Zb$, by saying that a bi-homogeneous element
of bi-degree $(i,j)$ has $\ZM$-degree $i-j$. If $c \in \CC$, then the ideal $\CG_c$ is not
bi-homogeneous (except if $c=0$) but is $\ZM$-homogeneous. Therefore,
$\Hb_c$ and $\Zb_c$ do not inherit a bi-grading, but they inherit
a $\ZM$-grading. Note finally that,
by definition, $\Pb=\CM[\CC] \otimes \CM[V]^W \otimes \CM[V^*]^W$ is clearly a bi-graded
subalgebra of $\Zb$.

The $\ZM$-grading on $\Zb$ and on $\Pb$ induces a $\CM^\times$-action on $\ZC$ and on $\PC$, and the morphism $\Upsilon \colon \ZC \to \PC$ is $\CM^\times$-equivariant.

\bigskip

\subsubsection{Poisson structure} \label{sec:poisson_bracket}
Let $t \in \CM$. One can define a deformation $\Hb_{t}$ of $\Hb$ as follows:
$\Hb_{t}$ is the quotient
of the algebra $\CM[\CC] \otimes (\Trm(V\times V^*)\rtimes W)$
by the relations
\equat\label{eq:rels-1}
\begin{cases}
[x,x']=[y,y']=0,\\
[y,x]=t \la y,x \ra + \DS{\sum_{s \in \REF(W)} (\e(s)-1)C_s
\frac{\langle y,\a_s \rangle \cdot \langle \a_s^\ve,x\rangle}{\la \a_s^\ve,\a_s \ra}\, s},
\end{cases}
\endequat
for all $x$, $x'\in V^*$, $y$, $y'\in V$. It is well-known~\cite{EG}
that the PBW-decomposition still holds so
that the family $(\Hb_{t})_{t \in \CM}$ is a flat deformation of $\Hb=\Hb_{0}$.
This allows one to define a Poisson bracket $\{\ ,\ \}$ on $\Zb$ as follows.
If $z_1$, $z_2 \in \Zb$,
we denote by $z_1^{t}$, $z_2^t$ the corresponding element of $\Hb_{t}$ through the
PBW-decomposition and we set
\equat \label{poisson_bracket}
\{z_1,z_2\} = \lim_{t \to 0} \frac{[z_1^t,z_2^t]}{t} .
\endequat
The Poisson bracket is $\CM[\CC]$-linear.

\bigskip

\subsubsection{Euler element}
Let $(y_1,\dots,y_n)$ be a basis of $V$ and let $(x_1,\dots,x_n)$
denote its dual basis. As in~\cite[\S{4.1}]{calogero}, we set
$$\euler  = ~\sum_{j=1}^n x_j y_j + \sum_{s \in \REF(W)} \e(s) C_s\, s
=~\sum_{j=1}^n x_j y_j +
\sum_{H\in\AC}\sum_{j=0}^{e_H-1}e_H\ K_{H,j}\varepsilon_{H,j}.$$
Recall that $\euler$ does not depend on the choice of the basis of $V$.
Also
\equat\label{eq:euler centre}
\euler \in \Zb, \quad \Frac(\Zb)=\Frac(\Pb)[\euler] \;,
\endequat
and
\equat\label{eq:euler poisson}
\{\euler,z\}=d z
\endequat
if $z \in \Zb$ is $\ZM$-homogeneous of degree $d$
(see for instance~\cite[Prop.~3.3.3]{calogero}).

\def\kb{\CM}

\bigskip
\subsection{Symplectic singularities and symplectic leaves} \label{sec:symp_sing}

The following concept has been introduced by Beauville \cite{beauville}: a normal variety $X$ is said to have {\it symplectic singularities} if its  smooth locus $X_{\mathrm{sm}}$ carries a symplectic form $\omega$ and for some (thus any) resolution of singularities $\pi : Y \to X$ the pull-back $\pi^*(\omega)$ defined on $\pi^{-1}(X_{\mathrm{sm}}) \subseteq Y$ extends to a (possibly degenerate) 2-form on all of $Y$.

For any $c \in \CC$ the Calogero--Moser space $\ZC_c$ is a variety with symplectic singularities
by~\cite{gordon icra}. Moreover, $\ZC_c$ admits a stratification into smooth symplectic subvarieties called the {\it symplectic leaves} of $\ZC_c$, see \cite{BrGo}. If $\mG$ denotes the maximal ideal in $\Zb_c$ corresponding to a point $p \in \ZC_c$, then the symplectic leaf $\mathcal{L}_c(p)$ containing $p$ is the subvariety defined by the largest ideal $I$ contained in $\mG$ which is Poisson, i.e. $\{Z_c,I\} \subset I$.

\def\Lie{{{\LG\iG\eG}}}

The zero-dimensional symplectic leaves are also called {\it cuspidal} points. A cuspidal point $p \in \ZC_c$ corresponds to a maximal ideal $\mG$ in $\Zb_c$ such that $\{\mG,\mG\} \subset \mG$. It follows that the cotangent space $\mG/\mG^2$ in $p$ inherits  a structure of a Lie algebra from the Poisson bracket. We will denote this Lie algebra by $\Lie_c(p)$.

\bigskip
\subsection{CHAMP}

The rational Cherednik algebra is implemented in CHAMP. This means concretely that it is realized as an algebra with a basis given by the PBW-decomposition \eqref{eq:pbw} and that products of elements will be rewritten in the basis (the algorithm is described in \cite{champ}). Instead of working over the complex numbers it is (mostly) sufficient to work over the defining field of the complex reflection group $W$, which is a number field. CHAMP comes along with an extensive documentation and we can just give a minimal impression with the simple examples included in this paper. We note that because all actions in MAGMA are from the {\it right}, it was conceptually more consistent to implement in CHAMP the {\it opposite} algebra of the rational Cherednik algebra as defined here—this is only relevant when translating elements between paper and CHAMP.

\bigskip

\begin{exemple} \label{ex:champ1}
  We verify that the Euler element in the generic rational Cherednik algebra for the Weyl group of type $B_2$ is indeed central (this simply works by checking whether the element commutes with the generators $x_j,y_j$ and with a set of generators of the group).
\begin{lstlisting}
> W := ComplexReflectionGroup(2,1,2);
> H := RationalCherednikAlgebra(W, 0 : Type:="BR");
> eu := EulerElement(H);
> IsCentral(eu);
true
\end{lstlisting}
\end{exemple}


\bigskip
\section{Generators and presentation of the center}

We describe algorithms for computing a minimal system of algebra generators and a presentation of the center $\Zb$ of~$\Hb$. Since $\Zb_c=\Zb/\CG_c \Zb$, we automatically obtain in this way also a (minimal) system of generators and a presentation of $\Zb_c$ for any $c \in \CC$. We use these algorithms in the sequel to obtain new results about the Calogero--Moser space for $G_4$ (Section \ref{sec:cm_g4}) and to compute families and hyperplanes (Section~\ref{sec:families}) for many exceptional complex reflection groups.

\bigskip
\subsection{The truncation map}

The backbone of our algorithms is the computation of the inverse of the truncation map: by the PBW-decomposition there is for any $h \in \Hb$
a unique family of elements $(h_w)_{w \in W}$ of
$\kb[\CC] \otimes \kb[V] \otimes \kb[V^*]$ such that
$$h=\sum_{w \in W} h_w w.$$
We define the $\kb[\CC]$-linear map $\troncation : \Hb \to \kb[\CC] \otimes \kb[V] \otimes \kb[V^*]$
by
$$\troncation(h)=h_1.$$
It is easily seen that
the map $\troncation : \Hb \to \kb[\CC] \otimes \kb[V] \otimes \kb[V^*]$ is $W$-equivariant and
$\Pb$-linear.
Recall that $\CG_0$ denotes the ideal of $\kb[\CC]$ generated by $(C_s)_{s \in \refw}$.

\bigskip

\begin{lem}\label{lem:0}
The restriction of the map $\troncation$ to $\Zb$ induces an isomorphism of bi-graded $\kb[\CC]$-modules
$$\troncation : \Zb \longiso \kb[\CC \times V \times V^*]^W$$
which satisfies
$$\troncation(z)-z \in \CG_0 \Hb.$$
In other words, if $f \in \kb[\CC \times V \times V^*]^W \subset \kb[\CC] \otimes \kb[V] \otimes \kb[V^*]$,
then there exists a unique element $z=\sum_{w \in W} z_w w$ of $\Zb$ (where
$z_w \in \kb[\CC] \otimes \kb[V] \otimes \kb[V^*]$) such that
$$\begin{cases}
z_1 = f,\\
z_w \equiv 0 \mod \CG_0 \otimes \kb[V] \otimes \kb[V^*] & \text{if $w \neq 1$}.
\end{cases}$$
\end{lem}

\bigskip

\begin{proof}
The map $\troncation : \Zb \longiso \kb[\CC \times V \times V^*]^W$ is a $\Pb$-linear
map between two free $\Pb$-modules of rank $|W|$. Since $\Zb/\CG_0\Zb=\Zb_0$, it is surjective
modulo $\CG_0$. By the graded Nakayama lemma (see for instance~\cite[Lemma~2.24]{broue})
applied to graded modules over the graded ring $\kb[\CC]$, this shows that it is surjective.
As any surjective map between finitely generated free modules of the same rank is an isomorphism,
the proof is complete.
\end{proof}

\bigskip

Note for future reference the following consequence:

\bigskip

\begin{coro}\label{coro:puissance de c0}
Let $z \in \Zb$, $z \neq 0$ and let $m$ be such that
$z \in \CG_0^m Z$. Write $z=\sum_{w \in W} z_w w$, with
$z_w \in \kb[\CC] \otimes \kb[V] \otimes \kb[V^*]$.
Then $z_w \in \CG_0^{m+1} Z$ for all $w \in W \setminus \{1\}$.
\end{coro}

\bigskip

\begin{proof}
The map $\troncation$ is $\kb[\CC]$-linear, so it is sufficient
to prove the result for $m=0$. But this follows from the fact
that $\Zb/\CG_0\Zb = \Zb_0$.
\end{proof}

\bigskip

We shall now describe how to compute the inverse map
$$\troncation^{-1} : \kb[\CC \times V \times V^*]^W \to \Zb.$$
If $s \in \REF(W)$ and $\ph \in \kb[V]$, we set
$$\D_s(\ph)=\frac{\ph-s(\ph)}{\a_s} \in \kb[V].$$
If $f \in \kb[\CC] \otimes \kb[V] \otimes \kb[V^*] \subset \Hb$, then
\equat\label{eq:yf}
[y,f]=\sum_{s \in \REF(W)} \e(s) \langle y,\a_s \rangle C_s (\Id_{\kb[\CC]} \otimes \D_s \otimes s)(f) s.
\endequat
A proof of this relation can be found in \cite[\S{3.6}]{gordon}

\bigskip

We define
\equat
V^\reg = V \setminus \bigcup_{H \in \AC} H = \{ v \in V \mid wv \neq v \ \text{for all } w \in W \}.
\endequat
Let $y_{\reg} \in V^{\reg}$ and use the notation of Lemma~\ref{lem:0}. Then
$$[y_\reg,z]=0$$
and, by~(\ref{eq:yf}), we get
$$\sum_{w \in W} \Bigl(z_w (y_\reg-w(y_\reg)) w + \sum_{s \in \REF(W)} \e(s) \langle y_\reg,\a_s \rangle C_s
(\Id_{\kb[\CC]} \otimes \D_s \otimes s)(z_w) sw \Bigr)= 0.$$
Therefore,
\equat\label{eq:algo}
z_w (y_\reg - w(y_\reg)) = -\sum_{s \in \REF(W)} \e(s) \langle y_\reg,\a_s \rangle C_s
(\Id_{\kb[\CC]} \otimes \D_s \otimes s)(z_{s^{-1}w})
\endequat
for all $w \in W$.
If we denote by $z_w^{<r>}$ the reduction modulo $\CG_0^{r+1} \otimes \kb[V] \otimes \kb[V^*]$ of $z_w$,
the formula~(\ref{eq:algo}) implies (by reduction modulo $\CG_0^{r+1}$) that
\equat\label{eq:algo-rec}
z_w^{<r>} (y_\reg - w(y_\reg)) = -\sum_{s \in \REF(W)} \e(s) \langle y_\reg,\a_s \rangle C_s
(\Id_{\kb[\CC]} \otimes \D_s \otimes s)(z_{s^{-1}w}^{<r-1>})
\endequat
for all $w \in W$.
This formula yields the following algorithm for computing $\troncation^{-1}$.

\bigskip

\refstepcounter{theo}
\noindent{\bfit Algorithm~\thetheo.\label{algo:troncation}
Computation of ${\mathbf{Trunc^{-1}}}$}
\begin{leftbar}
Let $f \in \kb[\CC \times V \times V^*]^W$ be bi-homogeneous of bi-degree $(d,e)$ and set $\d=\min(d,e)$.
Let $z=\troncation^{-1}(f)$ and write $z=\sum_{w \in W} z_w w$ with $z_w \in \kb[\CC] \otimes \kb[V] \otimes \kb[V^*]$.
Then:
\begin{itemize}
\itemth{1} For all $w \in W$, we define inductively a sequence $(z_w^{[r]})_{0 \le r \le \d}$
of elements of $\kb[\CC] \otimes \kb[V] \otimes \kb[V^*]$ as follows:
$z_1^{[r]}=f$ for all $r$ and, if $w \neq 1$, then
$$\begin{cases}
z_w^{[0]} = 0,\\
z_w^{[r]} (y_\reg - w(y_\reg)) =
\DS{-\sum_{s \in \REF(W)} \e(s) \langle y_\reg,\a_s \rangle C_s
(\Id_{\kb[\CC]} \otimes \D_s \otimes s)(z_{s^{-1}w}^{[r-1]})}.
\end{cases}$$
Note that, at each step, $z_w^{[r]}$ is bi-homogeneous of bi-degree $(d,e)$.
\itemth{2} The formula~(\ref{eq:algo}) implies that the right-hand side of the above formula
is divisible by $y_\reg-w(y_\reg)$, which is non-zero. It also implies that
$z_w^{[r]} \equiv z_w \mod \CG_0^{r+1} \otimes \kb[V] \otimes \kb[V^*]$.
\itemth{3} Since $\d=\min(d,e)$, we get $z_w^{[\d]}=z_w$.
\end{itemize}
Therefore, $z=\sum_{w \in W} z_w^{[\d]} w$.
\end{leftbar}


\bigskip

\begin{exemple}
Continuing Example \ref{ex:champ1}, we compute $\troncation^{-1}(f)$ for $f = x_1y_1 + y_2 x_2$. The result is the Euler element. Recall that CHAMP implements the opposite of the rational Cherednik algebra and therefore the variables in $f$ are reversed.

\begin{lstlisting}
> R<y1,y2,x1,x2> := PolynomialRing(Rationals(),4);
> f := y1*x1 + y2*x2;
> fH := TruncationInverse(H, f);
> fH eq EulerElement(H);
true
\end{lstlisting}
\end{exemple}

%


\bigskip
\subsection{System of generators}

Now that we can compute $\troncation^{-1}$ it is easy to find a system of generators of $\Zb$.

\bigskip

\begin{prop}\label{prop:generateurs}
Let $(z_i^{(0)})_{i \in I}$ be a family of bi-homogeneous generators of the $\kb$-algebra
$\kb[V \times V^*]^W$. Let $z_i = \troncation^{-1}(z_i^{(0)})$ (for $i \in I$).
Then $(z_i)_{i \in I}$ is a family of bi-homogeneous generators of the $\kb[\CC]$-algebra $\Zb$.
\end{prop}

\bigskip

\def\red{{\mathrm{red}}}
\begin{proof}
Let $\red_0 : \Zb \to \kb[V \times V^*]^W$ denote the reduction modulo $\CG_0$. Then, by construction,
$\red_0(z_i)=z_i^{(0)}$ and so the result follows from the graded Nakayama Lemma.
\end{proof}

\bigskip

\refstepcounter{theo}
\noindent{\bfit Algorithm~\thetheo.\label{algo:generateurs}  Minimal system of generators}
\begin{leftbar}
A minimal system of generators of $\Zb$ is obtained as follows:
\begin{itemize}
\itemth{1} Compute a minimal system of bi-homogeneous generators $(z_i^{(0)})_{i \in I}$ of the $\kb$-algebra $\kb[V \times V^*]^W$ (this is ``just'' computational invariant theory and there are algorithms for this in MAGMA already, specifically the algorithm \cite{king} by King which is implemented via the command \texttt{FundamentalInvariants}).
\itemth{2} Use Algorithm~\ref{algo:troncation} to compute $z_i=\troncation^{-1}(z_i^{(0)})$.
\end{itemize}
Then $(z_i)_{i \in I}$ is a minimal system of bi-homogeneous generators of $\Zb$.
\end{leftbar}

\bigskip

\begin{exemple} \label{ex:champ2}
We continue Example \ref{ex:champ1} and compute a minimal system of generators of~$\Zb$. In CHAMP, a minimal system of bi-homogeneous generators of $\kb[V \times V^*]^W$ is computed via the command {\tt SymplecticDoublingFundamentalInvariants} and the computation of the generators of $\Zb$ via Algorithm \ref{algo:generateurs} proceeds and stores the elements in the same order. In this example there are 8 generators, their bi-degrees are as listed below.
\begin{lstlisting}
> Zgen := CenterGenerators(H);
> SymplecticDoublingFundamentalInvariants(W)[2];
y1*x1 + y2*x2
> Zgen[2] eq EulerElement(H);
true
> [ Bidegree(f) : f in SymplecticDoublingFundamentalInvariants(W) ];
[ <0, 2>, <1, 1>, <2, 0>, <0, 4>, <1, 3>, <2, 2>, <3, 1>, <4, 0> ]
\end{lstlisting}
\end{exemple}

\bigskip

\begin{rema}
We discovered a bug in the computation of fundamental invariants in MAGMA that resulted in the return value not necessarily being a {\it minimal} system of generators (which contradicts the definition of a system of fundamental invariants and also contradicts the algorithm in \cite{king}). After reporting this bug, we were informed by MAGMA developer Allan Steel that this bug existed since V2.22 (May 2016). It was fixed after our report in V2.26-9 (October 2021).\finl
\end{rema}

\bigskip
\subsection{Presentation}
We now come to the computation of a presentation of $\Zb$, i.e. to the computation of the relations between a (minimal) system of generators of $\Zb$.



We fix a minimal
system of bi-homogeneous generators $(z_i^{(0)})_{i \in I}$ of $\kb[V \times V^*]^W$ and we
set $z_i=\troncation^{-1}(z_i^{(0)})$ for $i \in I$. We fix a family of indeterminates $(\zb_i)_{i \in I}$
and we denote by $\pi_0 : \kb[(\zb_i)_{i \in I}] \longto \kb[V \times V^*]^W$ (resp.
$\pi : \kb[\CC] \otimes \kb[(\zb_i)_{i \in I}] \longto \Zb$) the unique morphism of $\kb$-algebras
(resp. $\kb[\CC]$-algebras) such that $\pi_0(\zb_i)=z_i^{(0)}$ (resp.
$\pi(\zb_i)=z_i$). We endow $\kb[(\zb_i)_{i \in I}]$ with the $(\NM \times \NM)$-grading such that
$\zb_i$ has the same bi-degree as $z_i$ (or $z_i^{(0)}$), so that $\pi_0$ and $\pi$ are bi-graded
morphisms.

%


By definition, $\pi_0$ is surjective and by Proposition~\ref{prop:generateurs},
$\pi$ is surjective. Our aim is to compute a (minimal) set
of generators of $\Ker(\pi)$. We first explain how to compute preimages under the surjective
morphism $\pi$.

\newpage

\refstepcounter{theo}
\noindent{\bfit Algorithm~\thetheo.\label{algo:preimages}  Preimages under $\pi$}
\begin{leftbar}
Let $z \in \Zb$ be bi-homogeneous of bi-degree $(d,e)$ and let $\d=\min(d,e)$.
A bi-homogeneous element of bi-degree $(d,e)$ in $\pi^{-1}(z)$ is computed
by induction on $\d$ as follows:
\begin{itemize}
\itemth{1} Let $z_0$ denote the image of $z$ in $\kb[V \times V^*]^W \simeq Z/\CG_0 Z$.
Let $F_0 \in \kb[(\zb_i)_{i \in I}]$ be bi-homogeneous such that $\pi_0(F_0)=z_0$
(this is ``just'' computational commutative algebra and there are algorithms for this in MAGMA already, specifically one can compute the relations between a system of generators of $\kb[V \times V^*]^W$ given by primary invariants and irreducible secondary invariants via the command \texttt{Relations} and can translate this to relations between the fundamental invariants using the
 function {\tt HomogeneousModuleTest}).
\itemth{2} Then $\pi(F) \equiv z \mod \CG_0 Z$. Write $\pi(F)-z=\sum_{s \in \refw} C_s h_s$,
with $h_s \in \Zb$, bi-homogeneous of bi-degree $(d-1,e-1)$. By induction, we can find a family
$(F_s)_{s \in \refw}$ of bi-homogeneous elements of $\kb[(\zb_i)_{i \in I}]$
of bi-degree $(d-1,e-1)$ such that $\pi(F_s)=h_s$.
\end{itemize}
Then $F + \sum_{s \in \refw} C_s F_s$ is bi-homogeneous of bi-degree $(d,e)$
and its image under~$\pi$ is equal to $z$.
\end{leftbar}

\bigskip

Using Algorithm~\ref{algo:preimages}, one can lift relations between
the $z_i^{(0)}$'s in $\kb[V \times V^*]^W$ to relations
between the $z_i$'s in $\Zb$:

\bigskip

\refstepcounter{theo}
\noindent{\bfit Algorithm~\thetheo.\label{algo:relations}  Relations}
\begin{leftbar}
Let $\r^{(0)} \in \Ker(\pi_0)$ be bi-homogeneous of bi-degree $(d,e)$ and let $\d=\min(d,e)$. Then:
\begin{itemize}
\itemth{1} $\pi(\r^{(0)}) \equiv 0 \mod \CG_0 Z$, so we can write $\pi(\r^{(0)})=\sum_{s \in \refw} C_s h_s$
where $h_s \in \Zb$ is bi-homogeneous of bi-degree $(d-1,e-1)$.
\itemth{2} By Algorithm~\ref{algo:preimages}, we can find a family
$(F_s)_{s \in \refw}$ of bi-homogeneous elements of $\kb[(\zb_i)_{i \in I}]$
of bi-degree $(d-1,e-1)$ such that $\pi(F_s)=h_s$.
\itemth{3} Set $\r=\r^{(0)}-\sum_{s \in \refw} C_s F_s$.
\end{itemize}
Then $\pi(\r)=0$ by construction, $\r$ is bi-homogeneous of bi-degree $(d,e)$ and
$\r \equiv \r^{(0)} \mod \CG_0 \otimes \kb[(\zb_i)_{i \in I}]$.
\end{leftbar}

\bigskip

The next theorem shows how one can use
Algorithms~\ref{algo:troncation},~\ref{algo:generateurs},~\ref{algo:preimages}, 
and~\ref{algo:relations} to obtain a presentation of the algebra
$\kb[\CC]$-algebra $\Zb$:

\bigskip
\bigskip

\begin{theo}\label{theo:presentation}
Let $(\r_j^{(0)})_{j \in J}$ be a family of bi-homogeneous
generators of the ideal $\Ker(\pi_0)$ and, for $j \in J$,
let $\r_j \in \Ker(\pi)$ be bi-homogeneous (of the same bi-degree as $\r_j^{(0)}$)
and such that $\r_j \equiv \r_j^{(0)} \mod \CG_0 \otimes \kb[(\zb_i)_{i \in I}]$
(such a $\r_j$ is produced by Algorithm~\ref{algo:relations}). Then $\Ker(\pi)$ is generated
by $(\r_j)_{j \in J}$.
\end{theo}

\bigskip
\bigskip

\begin{proof}
Let $A$ denote the quotient of $\kb[\CC] \otimes \kb[(\zb_i)_{i \in I}]$ by the ideal $\aG$
generated by the family $(\r_j)_{j \in J}$. The morphism $\pi$ induces a surjective morphism
of $\kb[\CC]$-algebras $\pit : A \surto \Zb$.
By construction, this morphism induces an isomorphism
$\pit_0 : A/\CG_0 A \longiso \Zb/\CG_0 \Zb \simeq \kb[V \times V^*]^W$.
Since $\Zb$ is a free $\kb[\CC]$-module, there exists a sub-$\kb[\CC]$-module $\Zb'$ of $A$ such that
$A = \Ker(\pit) \oplus \Zb'$. By reduction modulo $\CG_0$, we get that $\Ker(\pibt)/\CG_0\Ker(\pibt)=0$.
This forces $\Ker(\pit)=0$, thanks to the graded Nakayama Lemma.
\end{proof}

\bigskip

\begin{exemple}
We continue Example \ref{ex:champ1} and compute a presentation of~$\Zb$. In CHAMP, the command {\tt SymplecticDoublingInvariantRingPresentation} computes  the relations of the invariant ring and the computation of the relations for $\Zb$ via Algorithm \ref{algo:relations} proceeds and stores the relations in the same order. In this example, $\Zb$ is a quotient of a polynomial ring in 8 variables (corresponding to the 8 elements in the minimal system of generators that we computed previously) by an ideal generated by 9 relations.

\begin{lstlisting}[breaklines=true]
> Zpres := CenterPresentation(H);
> Universe(Zpres);
Polynomial ring of rank 8 over Polynomial ring of rank 2 over Rational Field
Order: Lexicographical
Variables: z1, z2, z3, z4, z5, z6, z7, z8
> Zpres;
[
3*z1^2*z3 - z1*z2^2 - z1*z6 + 2*C1^2*z1 + z2*z5 - 2*z3*z4,
-4*z1*z2*z3 + 2*z1*z7 + z2^3 + 2*z2*z6 - 4*C2^2*z2 - z3*z5,
2*z1*z8 + z2^2*z3 - 2*z2*z7 - z3*z6 + 2*C1^2*z3,
8*z1^3*z3 - 3*z1^2*z2^2 - 4*z1^2*z6 + (4*C1^2 + 8*C2^2)*z1^2 + 2*z1*z2*z5 - 8*z1*z3*z4 + 4*z1*z3*z8 + 2*z2^2*z3^2 + 2*z2^2*z4 - 4*z2*z3*z7 - 2*z3^2*z6 + 4*C1^2*z3^2 + 4*z4*z6 - 8*C2^2*z4 - z5^2,
-7*z1^2*z2*z3 + 6*z1^2*z7 + z1*z2^3 + 3*z1*z2*z6 + (2*C1^2 - 4*C2^2)*z1*z2 + 2*z2*z3*z4 - 4*z4*z7 - z5*z6 + 2*C1^2*z5,
8*z1^2*z3^2 - 8*z1^2*z8 - 10*z1*z2^2*z3 + 6*z1*z2*z7 + (8*C1^2 - 4*C2^2)*z1*z3 + 2*z2^4 + 3*z2^2*z6 + (-6*C1^2 - 8*C2^2)*z2^2 + z2*z3*z5 - 8*z3^2*z4 + 8*z4*z8 - 2*z5*z7 + (-4*C1^2 + 4*C2^2)*z6 + 8*C1^4 - 8*C1^2*C2^2,
-6*z1^2*z3^2 + 10*z1^2*z8 + 8*z1*z2^2*z3 - 8*z1*z2*z7 - z2^4 - 2*z2^2*z6 + (4*C1^2 + 4*C2^2)*z2^2 + 4*z3^2*z4 - 4*z4*z8 - z6^2 + 4*C1^2*z6 - 4*C1^4,
-4*z1*z2*z3^2 + 2*z1*z2*z8 + 4*z1*z3*z7 + 3*z2^3*z3 - 4*z2^2*z7 + z2*z3*z6 + (-2*C1^2 - 4*C2^2)*z2*z3 - 2*z3^2*z5 + 2*z5*z8 - 2*z6*z7 + 4*C1^2*z7,
-4*z1*z3^3 + 4*z1*z3*z8 - 2*z2^2*z3^2 - 2*z2^2*z8 + 8*z2*z3*z7 + 4*z3^2*z6 - 4*C2^2*z3^2 - 4*z6*z8 - 4*z7^2 + 8*C2^2*z8
]
\end{lstlisting}
\end{exemple}

\bigskip
\subsection{Poisson brackets on $\Zb$}
Algorithm \ref{algo:preimages} allows us to express Poisson brackets $\{ u,v \}$ of elements $u,v \in \Zb$ in terms of the generators $z_i$ so that we can algorithmically work with $\Zb$ as a Poisson algebra.

\bigskip

\refstepcounter{theo}
\noindent{\bfit Algorithm~\thetheo.\label{algo:brackets} Poisson brackets on $\Zb$}
\begin{leftbar}
Let $u,v \in \Zb$.
\begin{itemize}
\itemth{1} Compute the commutator $[u^t,v^t]$ in $\Hb_t$, see Section \ref{sec:poisson_bracket}, and project to $\Hb$ by specializing $t=0$. Denote this element by $w$. Note that $w \in \Zb$.

\itemth{2} Compute $w' = \pi^{-1}(z_3)$ using Algorithm \ref{algo:preimages}
\end{itemize}
Then $w'$ is an expression of $\{ u,v \}$ in terms of the generators $(z_i)_{i \in I}$ of $\Zb$.
\end{leftbar}

\bigskip

\begin{exemple}
We continue Example \ref{ex:champ1} and compute the Poisson brackets $\{z_i,z_j\}$ of all generators $z_i,z_j$ of $\Zb$. This is done in CHAMP with the command \texttt{PoissonMatrix}. We note that because of the opposite algebra convention, the Poisson brackets in CHAMP are the negatives of what we have on paper. Since $z_2 = \euler$ by Example \ref{ex:champ2}, the second row of the matrix gives the (negatives) of the Poisson brackets $\{z_i,\euler\}$. Recall that $\{z_i,\euler\} = \deg_\ZM z_i$ by \eqref{eq:euler poisson}.
\begin{lstlisting}
> pmatZ := PoissonMatrix(H);
> pmatZ[2];
( 2*z1     0 -2*z3  4*z4  2*z5     0 -2*z7 -4*z8)
\end{lstlisting}
\end{exemple}


\bigskip
\subsection{Summary of what is computable} \label{center_comp_overview}

Our Algorithm \ref{algo:troncation} for computing the inverse of the truncation map performs very well in practice. In particular, we were able to compute a complete minimal system of generators of $\Zb$ for all exceptional complex reflection groups except $G_{16} - G_{22}$ and $G_{27} - G_{37}$. All these elements are stored in the database of CHAMP and can be loaded quickly. To give an impression of the performance---but also of the complexity involved in these computations---we note that we were able to compute, e.g., the inverse image of a generator $z_i^{(0)}$ of $(\NM \times \NM)$-degree $(6,6)$ for the Weyl group of type $F_4$ (equal to the exceptional complex reflection group $G_{28}$) which has 6,010 monomials and the total number of monomials in the various coefficients of the inverse image $z_i$ in $\Zb$ is 569,936. This single computation took about 6 hours.

The higher the degree of the element of which we want to compute the inverse image, the more iterations Algorithm \ref{algo:troncation} will take---and each iteration will be more complicated than the previous. This is one of the reasons why we were not able to compute a complete system of generators of $\Zb$ for all the exceptional complex reflection groups. Another problem arises even before we can come to the iteration part of the algorithm: the invariant theory of the action of $W$ on $V \times V^*$ is very complicated and for the higher rank cases we were not even able to compute a system of generators of $\kb[V \times V^*]^W$. The first in the series of exceptional complex reflection groups where this was not possible anymore is the group $G_{30}$, which is also the Coxeter group of type $H_4$ and therefore very interesting.

In our algorithm for computing families and hyperplanes that we present in Section \ref{sec:families} we do not need a complete (minimal) system of generators of $\Zb$ but only the sub-system consisting of generators of $\mathbb{Z}$-degree 0. We were able to compute this also for $G_{27}$ and $G_{28}$. Unfortunately, for the interesting case $G_{30}$ we are so far not able to do this because we do not even know how many generators of $\mathbb{Z}$-degree 0 there are in the invariant ring. We think this case poses a very interesting problem---both theoretically and computationally.

Finally, we note that computing a presentation of $\Zb$ via Algorithm \ref{algo:relations} is computationally very challenging and works only in the smallest cases: dihedral groups up to order 16 and the exceptional group $G_4$. Nonetheless, even these few cases led to important insight and results, see Section \ref{sec:dihedral_groups} and Section \ref{sec:cm_g4}. In Table \ref{timings_dihedral} we provide an overview of the performance of the algorithms in the case of small dihedral groups.

\begin{table}[htbp]
\centering
\begin{tabular}{|c||c|c|c|c||c|c|c|}
  \hline
  $d$ & A & B & C & D & E & F & G \\
  \hline \hline
  3 & 0.02s & 7 & 1 & 0.02s & 0.03s & 5 & 0.05s \\ \hline
  4 & 0.02s & 8 & 2 & 0.02s & 0.04s & 9 & 0.42s \\ \hline
  5 & 0.03s & 9 & 2 & 0.15s & 0.09s & 14 & 4.33s \\ \hline
  6 & 0.01s & 10 & 3 & 0.12s & 0.18s & 20 & 29.9s \\ \hline
  7 & 0.07s & 11 & 3 & 1.38s & 0.39s & 27 & 211s \\ \hline
  8 & 0.06s & 12 & 4 & 4.00s & 15.5s & 35 & 9222s \\ \hline
\end{tabular}
\label{timings_dihedral}
\bigskip
\caption{Overview of the performance of the algorithms in the case of dihedral groups of order~$2d$ for $3 \le d \le 8$. The meaning of the columns is as follows. A:~computation time of $(z_i^{(0)})_{i \in I}$. B:~cardinality of $I$. C:~maximum of $\delta$ among the~$z_i^{(0)}$. D:~computation time of $(z_i)_{i \in I}$ via Algorithm \ref{algo:generateurs}. E:~computation time of the relations between the $z_i^{(0)}$. F:~number of relations. G:~computation time of the relations between the $z_i$ via Algorithm \ref{algo:relations}.}
\end{table}

\bigskip
\subsection{Case of dihedral groups} \label{sec:dihedral_groups}

We mention some new results about the Calogero--Moser space for dihedral groups that grew out of the fact that we were able to compute an explicit presentation of~$\Zb$ in small cases (see Table \ref{timings_dihedral}).

Assume that $\dim_\CM V = 2$ and
that $W$ is the dihedral group
of order $2d$, with $d \ge 3$. Write $W=\la s,t \ra$ where $s$, $t$ are
Coxeter generators. Fix $c \in \CC$ and write $a=c_s$ and $b=c_t$.
We assume that $ab \neq 0$. Recall that $\ZC_c$ has dimension $4$ and that
$a=b$ whenever $d$ is odd. Note the following facts:
\begin{itemize}
\item[$\bullet$] If $d \ge 4$ and $a=b$, then $\ZC_c$ has a unique singular point $p$,
which is necessarily cuspidal~\cite[Tab.~5.2]{bonnafe diedral}.

\item[$\bullet$] If $d \ge 6$ is even and $a \neq b$, then $\ZC_c$ has a unique singular point $q$,
which is necessarily cuspidal~\cite[Tab.~5.2]{bonnafe diedral}.
\end{itemize}
Recall from Section \ref{sec:symp_sing} that we denote by $\Lie_c(p)$ the Lie algebra in a cuspidal point $p$ of $\ZC_c$. If $\gG$ is a simple Lie algebra, we denote by $\OC_\mini(\gG)$ the minimal nilpotent orbit of $\gG$. The following facts have been proved in~\cite[Prop.~8.4~and~8.8]{bonnafe diedral}
and~\cite{BBFJLS}, and rely on the fact that we were able to compute an explicit presentation of the center $\Zb$:

\bigskip

\begin{prop}\label{prop:b2g2}
With the above notation:
\begin{itemize}
\itemth{a} If $d=4$ and $a=b$, then $\Lie_c(p) \simeq \sG\lG_3(\CM)$ and the
symplectic singularity $(\ZC_c,p)$ is equivalent to $(\overline{\OC}_\mini(\sG\lG_3(\CM)),0)$.

\itemth{b} If $d=6$ and $a \neq b$, then $\Lie_c(p) \simeq \sG\pG_4(\CM)$ and the
symplectic singularity $(\ZC_c,q)$ is equivalent to $(\overline{\OC}_\mini(\sG\pG_4(\CM)),0)$.
\end{itemize}
\end{prop}

\bigskip

The main application of our algorithms in case of dihedral groups
is given by the next theorem~\cite{BBFJLS}.
It has been obtained through an explicit description of a presentation of $\Zb_c$
whenever $a=b$ for any $d \ge 3$ (see~\cite{bonnafe diedral 2}) which was obtained after computing
the cases $3 \le d \le 8$ with our algorithms and finding some
general patterns. So, even though this does not appear finally in the proof of the main
results in~\cite{bonnafe diedral 2} and~\cite{BBFJLS},
it is fair to say that these two papers owe their existence to
the algorithms developed here. Before stating the result,
let us introduce some notation: we set $\gG_d=\sG\lG_2(\CM) \oplus {\mathrm{Sym}}^d(\CM^2)$,
and we endow it with the Lie algebra structure such that ${\mathrm{Sym}}^d(\CM^2)$ is a
commutative ideal of $\gG_d$ and the adjoint action of $\sG\lG_2(\CM)$ on ${\mathrm{Sym}}^d(\CM^2)$
coincides with the natural action.

\bigskip

\begin{theo}\label{theo:new}
If $d \ge 5$ and $a=b$, then $\Lie_c(p) \simeq \gG_d$ and the isolated symplectic singularity
$(\ZC_c,p)$ has trivial local fundamental group.
\end{theo}

\bigskip

This gives a new family of examples of isolated symplectic singularities
with trivial local fundamental group, answering a question of Beauville~\cite[(4.3)]{beauville}.

\bigskip
\section{The Calogero--Moser space for $G_4$} \label{sec:cm_g4}

In this section, we use our algorithms to investigate Calogero--Moser spaces for the exceptional complex reflection group $G_4$: we give an explicit presentation, identify the symplectic singularity in the origin, and confirm a conjecture about symplectic leaves.

\subsection{The group and parameters} \label{sec:g4_group}


We assume that $V=\CM^2$ and that
$W$ is the subgroup of $\Gb\Lb_2(\CM)=\Gb\Lb_\CM(V)$ generated by
$$
s=\begin{pmatrix}
1 & 0 \\
0 & \zeta \\
\end{pmatrix}
\qquad\text{and}\qquad
t=\frac{1}{3}\begin{pmatrix}
2\zeta + 1 & 2(\zeta - 1) \\
\zeta - 1 &  \zeta + 2 \\
\end{pmatrix},$$
where $\z$ is a primitive third root of unity. Then $W$ is a primitive
complex reflection group of type $G_4$. Note that
\equat\label{eq:g4}
|W|=24,\qquad|\Zrm(W)|=2,\qquad\text{and}\qquad \Nrm_{\Gb\Lb_\CM(V)}(W)=W \cdot \CM^\times.
\endequat
We denote by $(y_1,y_2)$ the canonical
basis of $\CM^2$ and by $(x_1,x_2)$ its dual basis.
Let $H_s$ denote the reflecting hyperplane
of $s$ (note that $|\AC/W|=1$ and $e_{H_s}=3$). For simplification, we set $K_j=K_{H_s,j}$.
Recall that $K_0+K_1+K_2=0$. We also fix $c \in \CC$ and set $k=\kappa(c)$ for the corresponding parameter in $\KC$ (see Section \ref{sec:parameters}). For simplification, we set
$k_j=k_{H_s,j}$. Note that $k_0+k_1+k_2=0$.

\bigskip

\subsection{Presentation of $\Zb$}
The computation of a presentation of $\Zb$ in CHAMP takes about 2 minutes and is performed by the following commands:
\begin{lstlisting}
> W:=ComplexReflectionGroup(4);
> H:=RationalCherednikAlgebra(W,0 : Type:="BR-K");
> Zpres:=CenterPresentation(H);
\end{lstlisting}
The presentation is also stored in the database and can be retrieved in an instance. This returns a family $(z_j)_{1 \le j \le 8}$ of $8$ generators
of $\Zb_k$ and a family $(\EC_j')_{1 \le j \le 9}$ of $9$
equations. To describe this presentation explicitly, we rename and reorder the generators as follows:
$$(X_1,Y_1,X_2,Y_2,A,B,C,\euler)=(z_2,z_5,z_6,z_8,z_3,z_4,z_7,z_1).$$
Moreover, we replace the equations $(\EC_j')$ by the equations $(\EC_j)$ defined as follows:
\begin{align*}
&\EC_1 = \EC_1' , \quad \EC_2 = -\frac{1}{2} \EC_4' , \quad \EC_3 = 2 \EC_5' , \quad \EC_4 = -\frac{1}{2} \EC_2' , \quad \EC_5 = -2\EC_3' , \quad \EC_6 = -\EC_6' , \\
&\EC_7 = \frac{1}{4}( \EC_7' - 8 \EC_9') , \quad \EC_8 = - \EC_8' - \EC_6' , \quad \EC_9 = -10\EC_9' .
\end{align*}
We then obtain the following theorem.

\bigskip
\bigskip

\begin{theo}\label{prop:centre-g4}
We have:
\begin{itemize}
\itemth{a} $\CM[V]^W=\CM[X_1,X_2]$ and $\CM[V^*]^W=\CM[Y_1,Y_2]$,
so that
$$\Pb=\CM[X_1,Y_1,X_2,Y_2].$$

\itemth{b} The $\ZM$-degrees of the generators are given by
$$\deg(X_1,Y_1,X_2,Y_2,A,B,C,\euler) = (4,-4,6,-6,2,-2,0,0).$$

\itemth{c} The $\CM[\CC]$-algebra $\Zb$ admits the following presentation:
$$\begin{cases}
\text{\it Generators:}~X_1,Y_1,X_2,Y_2,A,B,C,\euler.\\
\text{\it Relations: see Table~\ref{table:rel-g4}.}
\end{cases}$$
\end{itemize}
\end{theo}

\begin{table}[htbp]
{\small $$\begin{cases}
2 X_1 Y_1 - 15 \euler^4 + 702(K_1^2 + K_1 K_2 + K_2^2) \euler^2 + 12 \euler C
+ 2592 K_1 K_2 (K_1 + K_2) \euler + A B = 0,\\
2 X_1 Y_2 + 3 Y_1 \euler A - 9 \euler^3 B + 378(K_1^2 + K_1 K_2 + K_2^2) \euler B + 4 B C = 0,\\
9 X_1 \euler^3 - 324(K_1^2 + K_1 K_2 + K_2^2) X_1 \euler - 8 X_1 C + 2 X_2 B - 3 \euler A^2 = 0, \\
3 X_1 \euler B + 2 X_2 Y_1 - 9 \euler^3 A + 378(K_1^2 + K_1 K_2 + K_2^2) \euler A + 4 A C = 0,\\
9 Y_1 \euler^3 - 324(K_1^2 + K_1 K_2 + K_2^2) Y_1 \euler - 8 Y_1 C + 2 Y_2 A - 3 \euler B^2 = 0, \\
2 X_1^2 B - 3 X_1 \euler^2 A + 144(K_1^2 + K_1 K_2 + K_2^2) X_1 A + 10 X_2 \euler^3
- 468(K_1^2 + K_1 K_2 + K_2^2) X_2 \euler \\
\hphantom{AAAA}
- 8 X_2 C - 1728 K_1 K_2 (K_1 + K_2) X_2 + 2 Y_1^2 A - 3 Y_1 \euler^2 B
+ 144(K_1^2 + K_1 K_2 + K_2^2) Y_1 B \\
\hphantom{AAAA}
+ 10 Y_2 \euler^3 - 468 (K_1^2 + K_1 K_2 + K_2^2) Y_2 \euler - 8 Y_2 C
- 1728 K_1 K_2 (K_1 + K_2) Y_2 - A^3 - B^3 = 0,\\
9 X_1 Y_1 \euler^2 + 2 X_2 Y_2 - 27 \euler^6 + 11664 K_1 K_2 (K_1 + K_2) \euler^3
+ 61236(K_1^2 + K_1 K_2 + K_2^2)^2 \euler^2 \\
\hphantom{AAAA}
+ 2160(K_1^2 + K_1 K_2 + K_2^2) \euler C + 16 C^2 = 0,\\
2 Y_1^2 A - 3 Y_1 \euler^2 B + 144(K_1^2 + K_1 K_2 + K_2^2) Y_1 B + 10 Y_2 \euler^3
- 468 (K_1^2 + K_1 K_2 + K_2^2) Y_2 \euler - 8 Y_2 C \\
\hphantom{AAAA}
- 1728 K_1 K_2 (K_1 + K_2) Y_2 - B^3 = 0,\\
60  X_1 Y_1 \euler^2 + 1944(K_1^2 + K_1 K_2 + K_2^2)  X_1 Y_1 + 5  X_1 B^2 + 10  X_2 Y_2
+ 5 Y_1 A^2 - 360 \euler^6 + 280 \euler^3 C \\
\hphantom{AAAA}
+ 97200 K_1 K_2 (K_1+K_2) \euler^3 + 798984(K_1^2 + K_1 K_2 + K_2^2)^2 \euler^2
+ 14544(K_1^2 + K_1 K_2 + K_2^2) \euler C \\
\hphantom{AAAA}
+ 1819584 K_1 K_2 (K_1 + K_2) (K_1^2 + K_1 K_2 + K_2^2) \euler
+ 1332 (K_1^2 + K_1 K_2 + K_2^2) A B \\
\hphantom{AAAA}
- 17280 K_1 K_2 (K_1 + K_2) C = 0.\\
\end{cases}$$}
\refstepcounter{theo}
\caption{Relations for $\Zb$ in type $G_4$.}
\label{table:rel-g4}
\end{table}

\subsection{The symplectic singularity in the origin}

Further to the dihedral group case mentioned in Section \ref{sec:dihedral_groups}, we are able to identify a symplectic singularity in the Calogero--Moser space for $G_4$.

\bigskip

\begin{theo}\label{eq:g4-sl3}
Assume the parameter is $(k_0,k_1,k_2)=(0,1,-1)$. Then the origin is the unique singular point of $\ZC_c$ and the symplectic singularity $(\ZC_c,0)$ is equivalent to $(\overline{\OC}_\mini(\sG\lG_3(\CM)),0)$.
\end{theo}

\bigskip

\begin{rema}\label{rem:param-g4}
By the Harish-Chandra theory of symplectic leaves~\cite{bellamy cuspidal}, the Calogero--Moser
space $\ZC_c$ may admit an isolated singularity only if $(k_0+k_1)(k_0+k_2)(k_1+k_2)=0$
and $(k_0,k_1,k_2) \neq (0,0,0)$.
Since the equations are symmetric in the $k_i$'s, we must focus on the case
where $k_1+k_2=0$ and, after multiplying by a scalar, the only possible case is the
one considered in the above Theorem.\finl
\end{rema}

\bigskip

\begin{proof}
  The proof is obtained via explicit computations. First, we specialize the presentation of $\Zb$ in the fixed parameter $c$ to obtain a presentation of $\Zb_c$ and then create $\ZC_c$ as a scheme:
  \begin{lstlisting}
  //Specialization morphism for the fixed parameter
  > phi := hom<BaseRing(H)->Rationals() | [1,-1] >;
  //Base ring of relations of Zc
  > A8ring<X1,Y1,X2,Y2,A,B,C,eu> := PolynomialRing(Rationals(), 8);
  //Variable change morphism
  > psi := hom<CenterSpace(H) -> A8ring | [eu, X1, A, B, Y1, X2, C, Y2]>;
  //Presentation of Zc
  > Zcpres := [ psi(ChangeRing(f, phi)) : f in CenterPresentation(H) ];
  > A8 := AffineSpace(A8ring);
  > Zc := Scheme(A8,Zcpres);
  \end{lstlisting}
  Next, we compute the singular locus of $\ZC_c$ as a (reduced) scheme:
  \begin{lstlisting}
  > Zcsing := SingularSubscheme(Zc);
  > Zcsing := ReducedSubscheme(Zcsing);
  > Zcsing := Scheme(A8,MinimalBasis(Zcsing));
  \end{lstlisting}
  Now, we show that $0 \in \ZC_c \subset \CM^8$ is the unique singular point of $\ZC_c$.
  \begin{lstlisting}
  > Zcsing;
  Scheme over Rational Field defined by
  eu,
  C,
  B,
  A,
  Y2,
  X2,
  Y1,
  X1
  \end{lstlisting}
  Note that even though we compute over $\mathbb{Q}$, the result shows that also over $\mathbb{C}$ the singular locus just consists of the origin. Next, we show that the projective tangent cone of $\ZC_c$ at $0$ is smooth (the computation takes about 2 minutes):
  \begin{lstlisting}
  > cone := TangentCone(Zc, Zc ! [0,0,0,0,0,0,0,0]);
  > projcone := Scheme(Proj(CoordinateRing(A8)),MinimalBasis(cone));
  > IsSingular(projcone);
  false
  \end{lstlisting}
  It thus follows from~\cite[Intro.]{beauville} that the symplectic singularity $(\ZC_c,0)$ is equivalent to $(\overline{\OC}_\mini(\gG),0)$ for some simple Lie algebra $\gG$. Now, $\gG$ is the tangent space of
  $\overline{\OC}_\mini(\gG)$ at $0$. Its dimension is equal to $8$
  by the following command:
  \begin{lstlisting}
  > Dimension(TangentSpace(Zc, Zc ! [0,0,0,0,0,0,0,0]));
  8
  \end{lstlisting}
  But $\sG\lG_3(\CM)$ is the only simple Lie algebra of dimension $8$, so we have proved the claim.
\end{proof}

\bigskip
\subsection{Confirming a conjecture about symplectic leaves} \label{leaf_closure_g4}

For the moment, let $W$ be an arbitrary complex reflection group. For a variety $\XC$ we denote by $\XC^\nor$ its normalization. If moreover $\XC$ is affine and irreducible, then $\XC^\nor$ is also affine and $\CM[\XC^\nor]$ is the integral closure of $\CM[\XC]$ in its fraction field  $\CM(\XC)$.

We fix an element of finite order $\t$ of $\Nrm_{\Gb\Lb_\CM(V)}(W)$ and a parameter $c \in \CC$ such that $\lexp{\t}{c}=c$. Then $\t$ also acts on the Calogero--Moser space $\ZC_c$ and we denote
by $\ZC_c^\t$ the closed subvariety of fixed points under the action of $\t$
(endowed with its reduced structure). As explained in~\cite{auto},
the variety $\ZC_c^\t$ inherits a partition into symplectic leaves
and we recall a conjecture about their structure~\cite[Conj.~B]{auto}:

\bigskip

\begin{conjecture}[Bonnaf\'e]\label{conj:feuilles}
Let $\LC$ be a symplectic leaf of $\ZC_c^\t$ and let $\LCov$ denote its closure in $\ZC_c^\t$.
Then $\overline{\LC}^\nor$
is isomorphic, as a Poisson variety endowed with a $\CM^\times$-action, to
some Calogero--Moser space associated with another pair $(V_\LC,W_\LC)$
and some parameter $k_\LC \in \aleph(W_\LC)$.
\end{conjecture}

\bigskip

In fact,~\cite[Conj.~B]{auto} is somewhat more precise, as it explains how
to recover the pair $(V_\LC,W_\LC)$. However, the parameter $k_\LC$
is quite a mystery. Note that, even if $\t=\Id_V$, this conjecture
is still unproved. Let us give some known cases:
\begin{itemize}
\item[$\bullet$] If $k=0$ (see~\cite[\S{4}]{auto}: note that $k_\LC=0$ in this case).

\item[$\bullet$] If $\ZC_c$ is smooth, then $\ZC_c^\t$
is also smooth and its symplectic leaves are its irreducible components. It has been
proved in~\cite[Theo.~1.3]{bonnafe maksimau} that, whenever $\t$ is the scalar multiplication
by a root of unity, then these irreducible components are isomorphic, as varieties
endowed with a $\CM^\times$-action, to Calogero--Moser spaces. This proves part of
Conjecture~\ref{conj:feuilles} (it remains to check that the
isomorphism constructed in~\cite[Theo.~1.3]{bonnafe maksimau}
respects the Poisson structure).
Note however that precise formulas are given for the value of
$k_\LC$ in~\cite[Theo.~4.21$($b$)$]{bonnafe maksimau}.

\item[$\bullet$] If $W$ is a Weyl group of type $B$ and $\t=\Id_V$, then Conjecture~B
is proved by Bellamy-Maksimau-Schedler in~\cite{be-sc}. As explained
in~\cite[Coro.~10.7]{auto}, we can easily deduce from their result that this
also implies Conjecture~B for Weyl groups of type $D$, with $\t$ being the identity
or a non-trivial involutive graph automorphism.
\end{itemize}

\bigskip

We will now add one more example to this list:

\bigskip

\begin{theo}\label{theo:g4}
Assume that $W$ is of type $G_4$. Then Conjecture~\ref{conj:feuilles} holds.
\end{theo}

\bigskip

\begin{proof}
We assume that $W$ is the group $G_4$ as in Section \ref{sec:g4_group}. Note that $W$ acts trivially on $\ZC_c$, so replacing $\t$ by $w\t$ for some $w \in W$,
does not change the fixed point subvariety. Therefore, since the normalizer of $W$
is $W \cdot \CM^\times$, we may, and we will, assume that $\t$ is the scalar multiplication
by a root of unity. Since $-\Id_V \in W$, we may also assume that the order of $\t$
is even: let us denote it by $2d$. Then
$$\CM[\ZC_c^\t]=\Zb_c/\sqrt{I_\t},$$
where $I_\t$ is the ideal generated by those of the generators $(X_1,Y_1,X_2,Y_2,A,B,C,\euler)$
whose $\ZM$-degree is not divisible by $2d$. In particular, if $d \ge 4$,
then $I_\t=\langle X_1,Y_1,X_2,Y_2,A,B\ra$ and so $\ZC_c^\t$ is zero-dimensional.
In this case, there is nothing to prove. This means that we may, and we will,
assume that $d \in \{1,2,3\}$. 

\medskip

If $d \in \{2,3\}$, then the computation of $\ZC_c^\t$ has been done in~\cite[\S{5}]{bonnafe maksimau}
and it has been checked in both cases that there is only one irreducible component
of dimension $2$ (the other being of dimension $0$) and it has been checked that
it is isomorphic, as a variety endowed with a $\CM^\times$-action, to the Calogero--Moser space
associated with the pair $(\CM,\mub_{2d})$, as expected from~\cite[Conj.~B]{auto}.
However, it has not been checked in~\cite[\S{5}]{bonnafe maksimau} that the isomorphism
respects the Poisson bracket, but this can be easily done in CHAMP. We will explain
how precisely in the more difficult situation where $d=1$, see below.

\medskip

Consequently, this means that we may assume that $d=1$, so that $\ZC_c^\t=\ZC_c$.
We will now prove that the symplectic leaves of $\ZC_c$
satisfy Conjecture~\ref{conj:feuilles}.

Let $\LC$ be a symplectic leaf of $\ZC_c$.
First, note that $\dim(\ZC_c)=4$ so that $\LC$ has
dimension $0$, $2$ or $4$. If $\dim(\LC)=4$, then $\LCov=\ZC_c$
and there is nothing to prove. If $\dim(\LC)=0$, then $\LC$
is a point and there is again nothing to prove.
So we may assume that $\dim(\LC)=2$. Note also that
we may assume that $k \neq 0$ since the case $k=0$
has been treated for any group.

\bigskip

\def\sing{{\mathrm{sing}}}

By the theory of Bellamy~\cite[Prop.~4.9]{bellamy cuspidal},
$\LC$ is parametrized by a
pair $(P,p)$ where $P$ is a parabolic subgroup of $W$ and $p$ is a cuspidal
point of $\ZC_c(V/V^P,P)$. Since $\dim(\LC)=2$,~\cite[Prop.~4.7]{bellamy cuspidal}
forces $\dim(V^P)=1$ and so we may assume that
$P=\langle s \rangle$. But, by~\cite[\S{18.5.A}]{calogero},
$\ZC_c(V/V^P,P)$ admits a cuspidal point if and only if
at least two of the $k_j$'s are equal. As the isomorphism class
of $\ZC_c$ depends on $k_0$, $k_1$, $k_2$ only up
to permutation, we may assume that $k_1=k_2$. Since $k \neq 0$
and $k_0+k_1+k_2=0$ and since we can rescale the
parameters, this means that we may assume that
$$(k_0,k_1,k_2)=(-2,1,1).$$
As in the proof of Theorem \ref{eq:g4-sl3} and continuing with these computations, we create $\ZC_c$ as a scheme in CHAMP as follows:
\begin{lstlisting}
> phi := hom<BaseRing(H)->Rationals() | [1,1] >;
> Zcpres := [ psi(ChangeRing(f, phi)) : f in CenterPresentation(H) ];
> Zc := Scheme(A8,Zcpres);
\end{lstlisting}
Next, we compute the reduced singular subscheme $\ZC_c^\sing$
of $\ZC_c$ and check that it is of dimension $2$ and irreducible:
\begin{lstlisting}
> Zcsing:=SingularSubscheme(Zc);
> Zcsing:=ReducedSubscheme(Zcsing);
> Dimension(Zcsing);
2
> IsIrreducible(Zcsing);
true
\end{lstlisting}
We conclude that $\ZC_c^\sing$ is irreducible {\it over $\QM$}
but not necessarily geometrically irreducible. However,
this also implies at least that it is of pure dimension $2$.
But its irreducible components (over $\CM$) are symplectic
leaves of dimension $2$: since $\ZC_c(V/V^P,P)$
has only one cuspidal point,~\cite[Prop.~4.9]{bellamy cuspidal} implies that
$\ZC_c$ has only one symplectic leaf of dimension $2$,
and this forces $\ZC_c^\sing$ to be geometrically irreducible. In other words,
$$\ZC_c^\sing=\overline{\LC}.$$
Since $\Nrm_W(P)/P \simeq \mub_2$ acting on $V^P \simeq \CM$
by scalar multiplication, we need to find
$k_\LC \in \aleph(\Nrm_W(P)/P)=\aleph(\mub_2)$ and a $\CM^\times$-equivariant
isomorphism of Poisson varieties
$$\overline{\LC}^\nor \longiso \ZC_{k_\LC}(\CM,\mub_2).$$

Looking at the equations of $\Zb_c$ in CHAMP
we see that among them we have the two
following ones:
$$X_1^3=X_2^2\qquad\text{and}\qquad Y_1^3=Y_2^2.$$
This means that the elements $X=X_2/X_1$ 
and $Y=Y_2/Y_1$ belong to the integral closure of
the ring $\CM[\overline{\LC}]$. Let $A=\CM[\ZC_c^\sing][X,Y]$.
We create the affine scheme $\LC_0$ with coordinate ring $A$ in CHAMP as follows:
\begin{lstlisting}
> A6<A,B,C,X,Y,E>:=AffineSpace(Rationals(),6);
> L0:=Scheme(A6,[Evaluate(f,[X^2,Y^2,X^3,Y^3,A,B,C,E]) :
  f in Basis(Ideal(Zcsing))]);
\end{lstlisting}
We see that $\LC_0$ is not reduced and define $\LC_1$ to be its reduced subscheme:
\begin{lstlisting}
> IsReduced(L0);
false
> L1:=ReducedSubscheme(L0);
> MinimalBasis(L1);
[
X*Y - E^2 + 6*E + 72,
B - Y*E + 12*Y,
B*X - C + 45/2*E^2 - 351/2*E - 648,
A - X*E + 12*X
]
\end{lstlisting}
We see that $A$, $B$ and $C$ can be expressed as polynomials
in $X$, $Y$ and $\euler$. Therefore,
$$\LC_1 \simeq \{(x,y,e) \in \CM^3~|~(e-6)(e+12)=xy\}.$$
As it is a smooth, hence normal, variety this shows that
$\LC_1 = \overline{\LC}^\nor$. In other words, setting $e'=e+3$ we obtain:
\equat\label{feuille-g4}
\overline{\LC}^\nor \simeq \{(x,y,e') \in \CM^3~|~(e'-9)(e'+9)=xy\}.
\endequat
Now, let $k_\LC \in \aleph(\CM,\mub_2)$ be defined by $(k_\LC)_0=-9$
and $(k_\LC)_1=9$.
Note also that $X$ and $Y$ have degree $2$ and $-2$ respectively.
Using the description
of Calogero--Moser spaces in rank $1$ given by~\cite[Theo.~18.2.4]{calogero},
we get that we have a $\CM^\times$-equivariant isomorphism of
varieties
\equat\label{eq:feuille-g4-mu2}
\overline{\LC}^\nor \longiso \ZC_{k_\LC}(\CM,\mub_2).
\endequat
It remains to show that this isomorphism respects the Poisson
bracket. For degree reasons, it respects the Poisson
bracket with the Euler element, see~\eqref{eq:euler poisson}.
So, it remains to compute the Poisson bracket $\{X,Y\}$
on both sides. On the right-hand side of~(\ref{eq:feuille-g4-mu2}),
we have
$$\{X,Y\}=-4 \euler_\# + 12,$$
where $\euler_\#$ denotes the image of $\euler$
under the isomorphism~\eqref{eq:feuille-g4-mu2}.
On the left-hand side of~(\ref{eq:feuille-g4-mu2}), we have
\begin{multline*}
\{X,Y\}=\Bigl\{\frac{X_2}{X_1},\frac{Y_2}{Y_1}\Bigr\}=
\frac{1}{X_1^2Y_1^2}(X_1Y_1\{X_2,Y_2\} \\
-X_1Y_2\{X_2,Y_1\}
-X_2Y_1\{X_1,Y_2\}+X_2Y_2\{X_1,Y_1\}).
\end{multline*}
Since $\CM[\ZC_c^\sing]$ is an integral domain,
it is sufficient to show that
$$X_1Y_1\{X_2,Y_2\}
-X_1Y_2\{X_2,Y_1\}
-X_2Y_1\{X_1,Y_2\}+X_2Y_2\{X_1,Y_1\}=X_1^2Y_1^2(-4 \euler + 12)\leqno{(?)}$$
in $\CM[\ZC_c^\sing]$.
The matrix whose entries are the Poisson brackets $\{u,v\}$ for $u$ and $v$ running over the generators $X_1$, $Y_1$, $X_2$, $Y_2$, $A$, $B$, $C$, $\euler$ of $\Zb_c$ is obtained in CHAMP using Algorithm \ref{algo:brackets} as follows (recall that Poisson brackets in CHAMP are the negatives of our Poisson brackets because of the opposite algebra convention):
\begin{lstlisting}
> pmatZ := PoissonMatrix(H);
> pmatZc := Matrix(A8ring, 8, 8, [ -psi(ChangeRing(f,phi)) :
  f in Eltseq(pmatZ) ]);
//Permute rows and columns to account for change in ordering of generators
> pmatZc := Permute(pmatZc, [2,5,6,8,3,4,7,1]);
\end{lstlisting}
Finally, we confirm that equation~(?) indeed holds:
\begin{lstlisting}
> lhs := CoordinateRing(Zcsing)!(X1*Y1*pmatZc[3,4] - X1*Y2*pmatZc[3,2]
  - X2*Y1*pmatZc[1,4] + X2*Y2*pmatZc[1,2]);
> rhs := CoordinateRing(Zcsing)!(X1^2*Y1^2*(-4*eu+12));
> lhs eq rhs;
true
\end{lstlisting}
This completes the proof.
\end{proof}

\bigskip
\section{Families and hyperplanes} \label{sec:families}

\medskip

In this section we describe algorithms for computing the Calogero--Moser families (which correspond to the $\mathbb{C}^\times$-fixed points of $\ZC_c$) and an associated hyperplane arrangement which plays an important role in geometry (Section \ref{sec:min_models}) and representation theory (Section \ref{sec:martino}). We furthermore give an algorithm for computing the cuspidal points (i.e.\ the zero-dimensional symplectic leaves) of the Calogero--Moser space.

\subsection{Restricted algebras, families, and hyperplanes}
Throughout, let $W$ be an arbitrary complex reflection group. Let $\mG_{0}$ be the maximal ideal in $\CM[V]^W \otimes \CM[V^*]^W$ corresponding to the origin of $V/W \times V^*/W$. The $\kb[\CC]$-algebra $\overline{\Hb} = \Hb/\mG_0 \Hb$ is called the {\it generic restricted rational Cherednik algebra}. We can specialize $\overline{\Hb}$ in $c \in \CC$ to obtain the $\kb$-algebra $\overline{\Hb}_c = \overline{\Hb}/\CG_c \overline{\Hb}$. Since $\Hb$ is a free $\Pb$-module of rank $|W|^3$, it follows that $\overline{\Hb}$ is a free $\kb[\CC]$-module of rank $|W|^3$ and that $\overline{\Hb}_c$ is a $\kb$-algebra of dimension $|W|^3$. The algebras $\overline{\Hb}_c$ were studied by Gordon~\cite{gordon}.

The PBW-decomposition of $\Hb$ induces a triangular decomposition of $\overline{\Hb}_c$ and this leads to a theory of standard modules similar to the theory for finite-dimensional complex semisimple Lie algebras. Specifically, for each $\lambda \in \Irr(W)$ there is an associated $\overline{\Hb}_c$-module $\Delta_c(\lambda)$ which has simple head $L_c(\lambda)$ and $\{ L_c(\lambda) \mid \lambda \in \Irr W \}$ is a complete set of representatives of isomorphism classes of simple $\overline{\Hb}_c$-modules. The details are given in Gordon~\cite{gordon}. Since $\overline{\Hb}_c$ is a finite-dimensional $\kb$-algebra, it has a block decomposition and each simple $\overline{\Hb}_c$-module belongs to exactly one block. Via the bijection $\lambda \mapsto L_c(\lambda)$ this block decomposition yields a partition of $\Irr W$ that we denote by $\mathrm{CM}_c$. The parts of this partition are called the {\it Calogero--Moser $c$-families}.

Let $\Upsilon_c \colon \ZC_c \to \PC_c$ be the fiber of the morphism $\Upsilon \colon \ZC \to \PC$ from \eqref{upsilon} in $c$, i.e. $\Upsilon_c$ is the morphism defined by the inclusion $\Pb_c \subseteq \Zb_c$. This morphism is $\CM^\times$-equivariant and by \cite{gordon} the blocks of $\overline{\Hb}_c$, and thus the Calogero--Moser $c$-families, are naturally in bijection with the $\CM^\times$-fixed points of $\ZC_c$. Denoting by $\ZC_c^{\CM^\times}$ the set of $\CM^\times$-fixed points of $\ZC_c$ we define
\equat \label{cm_locus}
N = \max_{c \in \CC} \left| \ZC_c^{\CM^\times} \right| \quad \text{and} \quad \CC_\calo = \left\lbrace c \in \CC \mid  \left| \ZC_c^{\CM^\times} \right| < N \right\rbrace .
\endequat
The set $\CC_\calo$ is known to be a union of (finitely many) hyperplanes (see below): these are the {\it Calogero--Moser hyperplanes}.

\bigskip

\begin{rema}
The original purpose of CHAMP, as presented in \cite{champ}, was to compute decomposition matrices of  standard modules for restricted rational Cherednik algebras—from which one can easily deduce the Calogero--Moser families. The database of CHAMP contains a wealth of such information for several exceptional complex reflection groups. We note that the algorithm for computing Calogero--Moser families we will present below (Algorithm \ref{algo:familles}) is much more efficient and allows to cover many more cases than the computation via decomposing standard modules (which, on the other hand, yields much more information).\finl
\end{rema}

\bigskip
\bigskip
\subsection{Families and hyperplanes via central characters} \label{sec:families_via_central_char}
The Calogero--Moser families and hyperplanes admit an explicit description via central characters that—combined with the fact that we can compute generators of $\Zb$ via Algorithm \ref{algo:generateurs}—immediately yields an algorithm for computing these invariants. We denote by
$$\Omeb^\Hb : \Hb \longto \kb[\CC] \otimes \kb W$$
the unique $\kb[\CC]$-linear map which, through the isomorphism~(\ref{eq:pbw}) of
vector spaces given by the PBW-decomposition, sends an element $f \in \kb[V]$ (resp. $f^* \in \kb[V^*]$)
to $f(0)$ (resp. $f^*(0)$), and which is the identity on $\kb[\CC]$ and $\kb W$.
This map is $W$-equivariant, so its restriction to $\Zb$ induces a map
$$\Omeb : \Zb \longto \kb[\CC] \otimes \Zrm(\kb W) \;,$$
which is a morphism of $\kb[\CC]$-algebras by \cite[Coro.~4.2.11]{calogero}.
For a prime ideal $\CG$ of $\kb[\CC]$ we denote by
$$
\Omeb^{\CG} : \Zb \to (\kb[\CC]/\CG) \otimes \Zrm(\kb W)
$$
the composition of $\Omeb$ with the scalar extension of the quotient map $\kb[\CC] \to \kb[\CC]/\CG$ to $\Zrm(\kb W)$. Note that $\Omeb^{(0)} = \Omeb$ for the zero ideal $(0)$ in $\kb[ \CC]$. For $c \in \CC$ we set
$$\Omeb^c = \Omeb^{\CG_c} : \Zb \to \Zrm(\kb W) .
$$
For $\chi \in \Irr(W)$ we denote by $\o_\chi : \Zrm(\kb W) \longto \kb$
the associated central character and also denote by
$\o_\chi : (\kb[\CC]/\CG) \otimes \Zrm(\kb W) \longto \kb[\CC]/\CG$
its scalar extension to $\kb[\CC]/\CG$. We then set
\equat
\Omeb_\chi^{\CG} = \o_\chi \circ \Omeb^{\CG} : \Zb \to \kb[\CC]/\CG \;.
\endequat
It follows from \cite[Theorems 1.1 and 1.5]{thiel-blocks} that two characters $\chi,\chi'$ are in the same Calogero--Moser $c$-family if and only if $\Omeb_\chi^c = \Omeb_{\chi'}^c$. We will therefore say more generally that for a prime ideal $\CG$ in $\kb [\CC]$ two characters $\chi,\chi'$ are in the same {\it Calogero--Moser $\CG$-family} if $\Omeb_\chi^{\CG} = \Omeb_{\chi'}^{\CG}$. We denote by $\CM_{\CG}$ the set of Calogero--Moser $\CG$-families. The Calogero--Moser $\CG$-families for $\CG=(0)$ are called the {\it generic Calogero--Moser families}. These are determined by the condition $\Omeb_\chi = \Omeb_{\chi'}$. It is clear that if $\CG' \subseteq \CG$ is an inclusion of prime ideals of $\kb[\CC]$, then the $\CG$-families are unions of $\CG'$-families, i.e. specializing yields coarser families. In particular, the number of $\CG'$-families is less than or equal to the number of $\CG$-families. Let
\equat
\widetilde{\CC_\calo} = \{ \CG \in \mathrm{Spec}(\kb[\CC]) \mid \Omeb_\chi \neq \Omeb_{\chi'} \text{ but } \Omeb_\chi^\CG = \Omeb_{\chi'}^{\CG} \text{ for some } \chi,\chi' \}
\endequat
be the locus of all $\CG$ where the $\CG$-families are not equal to (and thus coarser than) the generic families. It follows from \cite[Theorem 1.3]{thiel-blocks} that this locus is a closed subset of $\mathrm{Spec}(\kb[\CC])$ which is either empty or pure of codimension 1. It is obviously contained in
 $$
 \widetilde{\CC_{\mathrm{Eu}}} = \{ \CG \in \mathrm{Spec}(\kb[\CC]) \mid \Omeb_\chi(\euler) \neq \Omeb_{\chi'}(\euler) \text{ but } \Omeb_\chi^\CG(\euler) = \Omeb_{\chi'}^{\CG}(\euler) \text{ for some }  \chi,\chi' \} .
$$
This set is the union of the zero sets of the polynomials $\Omeb_\chi(\euler) - \Omeb_{\chi'}(\euler) \in \kb[\CC]$ for characters $\chi,\chi'$ with $\Omeb_\chi(\euler) \neq \Omeb_{\chi'}(\euler)$. Note that $\Omeb_\chi(\euler) \in \kb[\CC]$ is a linear polynomial. It thus follows that $\widetilde{\CC_{\mathrm{Eu}}}$ is a union of (finitely many) hyperplanes: we call them the {\it Euler hyperplanes}. In particular, $\widetilde{\CC_{\mathrm{Eu}}}$ is non-empty and thus pure of codimension 1. This shows that there is $c \in \CC$ such that the $c$-families are equal to the generic families and this implies that the set $\CC_\calo$ from \eqref{cm_locus} is equal to (the set of closed points of) $\widetilde{\CC_\calo}$.   It has been proven in~\cite[Coro.~7.8.3]{calogero} that
$\widetilde{\CC_\calo}$ is a union of Euler hyperplanes. In particular, $\widetilde{\CC_\calo}$ is a union of hyperplanes.

Given two partitions $\mathcal{P},\mathcal{P}'$ of $\Irr W$ we denote by $\mathcal{P} \wedge \mathcal{P}'$ the partition of $\Irr W$ whose parts are the unions of all parts of $\mathcal{P}$ and $\mathcal{P}'$ having non-empty intersection. Let $H_1,\ldots,H_r$ be the Calogero--Moser hyperplanes and for each $j$ let $\CG_j$ be the corresponding prime ideal in $\kb[\CC]$. It follows from \cite[Lemma 3.6]{thiel-blocks} that for $c \in \CC_\calo$ we have
\equat \label{CM_semicont}
\mathrm{CM}_c = \bigwedge_{ \substack{j=1,\ldots,r \\ c \in H_j}} \mathrm{CM}_{\CG_j} \;,
\endequat
i.e. the $c$-families are obtained by taking the meet of the $\CG_j$-families for all $i$ such that $c$ is contained in $H_i$.

Our discussion yields an immediate algorithm for computing the Calogero--Moser hyperplanes and the $c$-families for {\it all} parameters $c \in \CC$ at once. Before stating this we record an elementary fact that simplifies the computations.

\bigskip
\begin{lem}\label{lem:zplus}
Let $z \in \Zb$ be $\ZM$-homogeneous of non-zero degree. Then
$\Omeb(z)=0$.
\end{lem}

\bigskip

\begin{proof}
Indeed, the morphism $\Omeb : \Zb \to \kb[\CC] \otimes \Zrm(\kb W)$ respects
the $(\NM \times \NM)$-bigrading, so respects the $\ZM$-grading. But $\kb[\CC] \otimes \Zrm(\kb W)$ is concentrated in $\ZM$-degree $0$.
\end{proof}

\bigskip

\refstepcounter{theo}
\noindent{\bfit Algorithm~\thetheo.\label{algo:familles} Families and hyperplanes}
\begin{leftbar}
  \begin{itemize}
  \itemth{1} Let $(z_i^{(0)})_{i \in I}$ be a system of bi-homogeneous generators
  of the $\kb[\CC]$-algebra $\kb[V \times V^*]^W$. Let $I_0=\{i \in I~|~\deg_\ZM(z_i^{(0)})=0\}$. Compute the corresponding elements $(z_i)_{i \in I_0}$ of $\Zb$ via Algorithm~\ref{algo:generateurs}.

  \itemth{2} Compute the generic Calogero--Moser families: these are the fibers of the map
  $$\fonctio{\Irr(W)}{\kb^{I_0}}{\chi}{\bigl(\Omeb_\chi(z_i)\bigr)_{i \in I_0}.}$$

  \itemth{3} Compute the Calogero--Moser hyperplanes: for any two characters $\chi, \chi'$ with $\Omeb_\chi \neq \Omeb_{\chi'}$ compute the primary decomposition of the radical of the ideal in $\kb[\CC]$ generated by the elements $\Omeb_\chi(z_i) - \Omeb_{\chi'}(z_i)$ for $i \in I_0$. In MAGMA this can be done via the command {\tt RadicalDecomposition}, which is more efficient than first computing the radical and then computing its primary decomposition. The collection of the hyperplanes obtained in this way are the Calogero--Moser hyperplanes.

  \itemth{4} Compute the Calogero--Moser families on the hyperplanes: let $H_1,\ldots,H_r$ be the Calogero--Moser hyperplanes and for each $j$ let $\CG_j$ be the corresponding prime ideal in $\kb[\CC]$. The $\CG_j$-families are the fibers of the map
  $$\fonctio{\Irr(W)}{\kb^{I_0}}{\chi}{\bigl(\Omeb_\chi^{\CG_j}(z_i)\bigr)_{i \in I_0}.}$$

  \itemth{5} For arbitrary $c \in \CC$ the Calogero--Moser $c$-families are computed as follows. Since $\CC_\calo$ was computed in (3), we can test whether $c \in \CC_\calo$. If $c \notin \CC_{\calo}$, the $c$-families are equal to the generic families computed in (2). If $c \in \CC_\calo$, the $c$-families are given by
  $$
  \mathrm{CM}_c = \bigwedge_{ \substack{j=1,\ldots,r \\ c \in H_j}} \mathrm{CM}_{\CG_j} \;,
  $$
  where $\mathrm{CM}_{\CG_j}$ was computed in (4).
  \end{itemize}
\end{leftbar}

\bigskip
\begin{exemple}
We compute in CHAMP the Calogero--Moser hyperplanes and families for the Weyl group of type $B_2$.
\begin{lstlisting}
> W:=ComplexReflectionGroup(2,1,2);
> H:=RationalCherednikAlgebra(W,0 : Type:="BR-K");
> hyp := CalogeroMoserHyperplanes(H); hyp;
[
K2_1,
K1_1,
K1_1 + K2_1,
K1_1 - K2_1
]
// The families on the hyperplane K1_1 - K2_1
> CalogeroMoserFamilies(H)[hyp[4]];
{
{ 1, 2, 5 },
{ 3 },
{ 4 }
}
\end{lstlisting}
\end{exemple}

\bigskip
\bigskip
\subsection{Summary of what is computable} \label{sec:families_summary}
Remember from Section \ref{center_comp_overview} that we could compute $(z_i)_{i \in I_0}$ for all exceptional complex reflection groups except $G_{16} - G_{22}$ and $G_{29} - G_{37}$. In all these cases, we could successfully run Algorithm \ref{algo:familles} and we thus know the Calogero--Moser hyperplanes and the $c$-families for all $c \in \CC$. The results are contained in the database of CHAMP. It is particularly exciting that we were able to compute this for $G_{28}$ because it is a Weyl group (of type $F_4$). Table \ref{table:cm_hyperplanes} gives a brief summary of results about the Calogero--Moser hyperplane arrangement for exceptional complex reflection groups (we have published this table previously in \cite{BST-Hyperplanes}). In the following two subsections we mention  important results that follow from these computations.

\begin{table}[htbp]
\begin{displaymath}
{\footnotesize
\begin{array}{|r||c|c|c|}
\hline
\textrm{Group} & \#\text{Hyp} & \textrm{Poincar\'e polynomial} & \QM\Frm\Trm \\
\hline \hline
G_4   & 6 & (5t+1)(t+1) & 2   \\ \hline
G_5   & 33 & (116t^2 + 21t + 1)(11t + 1)(t+1) & 92  \\ \hline
G_6   & 16  & (8t + 1)(7t + 1)(t+1)& 12  \\ \hline
G_7   &61 & (98644t^4 + 18462t^3 + 1489t^2 + 60t + 1)(t+1) & 3296  \\ \hline
G_8   & 25 & (13t + 1)(11t + 1)(t+1) & 14  \\ \hline
G_9   & 54 &  (6499t^3 + 983t^2 + 53t + 1)(t+1) & 2  \\ \hline
G_{10}  & 111 & (1001586t^4 + 107662t^3 + 4913t^2 + 110t + 1)(t+1) & 15476  \\ \hline
G_{11}  &  196 & (383999826t^5 + 25688824t^4 + 857259t^3 + 17047t^2 + 195t + 1)(t+1) & 2851133  \\ \hline
G_{13}  & 6    & (5t+1)(t+1) & 3 \\ \hline
G_{14}  &  22  & (116t^2 + 21t + 1)(t+1) & 23  \\ \hline
G_{15}  &  65 & (13982t^3 + 1529t^2 + 32t + 1)(1+t) & 2596  \\ \hline
G_{20}  & 12   & (11t+1)(t+1) & 4  \\ \hline
G_{25}  &  12   & (11t+1)(t+1) & 4 \\ \hline
G_{26}  &  37  & (335t^2 + 36t + 1)(t+1) & 62  \\ \hline
G_{28}   &  8   & (7t+1)(t+1) & 4  \\ \hline
\end{array}
}
\end{displaymath}
\label{table:cm_hyperplanes}
\caption{Some data about the Calogero--Moser hyperplane arrangement for exceptional complex reflection groups: the number of hyperplanes, the Poincaré polynomial of the arrangement, and the number of $\mathbb{Q}$-factorial terminalizations of the associated symplectic singularity (see Section \ref{sec:min_models}). We note that the groups $G_{12}$, $G_{22}$ -- $G_{24}$, $G_{27}$, $G_{29}$ -- $G_{31}$, and $G_{33}$ -- $G_{37}$ have a single conjugacy class of reflections so that $\CC$ is 1-dimensional and $\CC_{\calo}$ is just the origin. Hence, only for the exceptional complex reflection groups $G_{16}$ -- $G_{19}$, $G_{21}$, and $G_{32}$ the Calogero--Moser hyperplanes are still unknown.}
\label{exc_table}
\end{table}

\bigskip
\subsection{Applications in birational geometry} \label{sec:min_models}
Let $X = (V \times V^*)/W$. This variety is normal. Recall that the vector space $V \times V^*$ carries a natural symplectic form and the action of $W$ on this space is symplectic. Since the symplectic group is contained in the corresponding special linear group, it follows from the Reid--Tai criterion that $X$ has canonical singularities (see, e.g., \cite[Theorem 3.21]{kollar}). An application of \cite[Corollary 1.4.3]{BCHM} thus shows that $X$ admits a {\it $\mathbb{Q}$-factorial terminalization}, i.e. a crepant projective birational morphism $\pi \colon Y \to X$ from a normal $\mathbb{Q}$-factorial variety $Y$ with terminal singularities (the projectivity of the morphism follows from the proof of \cite[Corollary 1.4.3]{BCHM}). Since $X$ has trivial canonical class by a result by Watanabe (see, e.g. \cite[Theorem 4.6.2]{benson}), the $\mathbb{Q}$-factorial terminalizations are the same as {\it relative minimal models} of $X$ (relative to a log-resolution of~$X$).

Now, choose such a $\mathbb{Q}$-factorial terminalization $\pi \colon Y \to X$. For more details on the following we refer to \cite{Namikawa-Birational}. The cohomology group $\mathcal{K}_{\mathbb{C}} = H^2(Y,\mathbb{C})$ can be identified with $\Pic(Y) \otimes_{\mathbb{Z}} \mathbb{C}$. Let $\Movable(\pi)$ be the cone in $\mathcal{K}_{\mathbb{R}} = \Pic(Y) \otimes_{\mathbb{Z}} \mathbb{R}$ formed by the $\pi$-movable line bundles. This cone decomposes into the ample cones of the various other $\mathbb{Q}$-factorial terminalizations of $X$, and the codimension-1 faces of each of these ample cones generate a hyperplane arrangement in $\mathcal{K}_{\mathbb{R}}$ decomposing $\Movable(\pi)$ into various chambers.

It is shown in \cite{bellamy counting} that $\mathcal{K}_{\mathbb{R}}$ can naturally be identified with the real vector space inside the complex parameter space $\KC$ spanned by the $K_{\O,j}$, i.e., $\KC_\RM=\sum_{(\O,j) \in \aleph(W)} \RM K_{\O,j}$. Recall from Section \ref{sec:families_via_central_char} that the Calogero--Moser hyperplanes are Euler hyperplanes. It is then clear that in the parameters $K_{\O,j}$ the Calogero--Moser hyperplanes have rational coefficients and so they define a hyperplane arrangement in $\KC_\RM$ as well. Using the Poisson-geometric description of the chamber decomposition of the movable cone given in \cite{Namikawa-Birational}, it was shown in \cite{BST-Hyperplanes} that these chambers coincide with orbits of the chambers of the Calogero--Moser hyperplane arrangement under the action of the {\it Namikawa Weyl group} on $\KC_{\RM}$. This group was shown in \cite{bellamy counting} to be simply the product of the full permutation groups on the indices $j$ of the variables $\KC_\RM=\sum_{(\O,j)}$ for fixed $\O$.

To summarize, when we denote by $\Erm_\calo(W)$ the number of chambers of the Calogero--Moser hyperplane arrangement inside $\KC_\RM$ and by $\QM\Frm\Trm(W)$ the number of $\mathbb{Q}$-factorial terminalizations of $(V \times V^*)/W$, then
\equat\label{eq:bellamy}
\QM\Frm\Trm(W)=\frac{\Erm_\calo(W)}{\DS{\prod_{\O \in \AC/W} (e_\O !)}}.
\endequat

Recall that the number of chambers of a real hyperplane arrangement is equal to the evaluation at $-1$ of its Poincaré polynomial \cite{zaslavsky}. Algorithms for computing the Poincaré polynomial are implemented in, e.g., the computer algebra system Sage \cite{Sage}. In all cases of the exceptional complex reflection groups where we could compute the Calogero--Moser hyperplane arrangement we were also able to compute the number $\QM\Frm\Trm(W)$ as well, see Table \ref{table:cm_hyperplanes}. We note that the computation of the Poincaré polynomial in case of $G_{11}$ in Sage took three weeks.

\bigskip

\begin{rema}
It would be very interesting to explicitly construct a $\mathbb{Q}$-factorial terminalization of $(V \times V^*)/W$. Among the exceptional groups, this has so far only been done for $G_4$, see \cite{lehn-sorger}, which is very special because it is the only exceptional group where the associated symplectic singularity admits a symplectic resolution, i.e. the $\mathbb{Q}$-factorial terminalizations are smooth. Recently, an explicit presentation of a $\mathbb{Q}$-factorial terminalization in the case of odd dihedral groups was constructed in \cite[Corollary 7.3.13]{SchmittThesis}. \finl
\end{rema}

\bigskip
\bigskip
\subsection{Applications in representation theory} \label{sec:martino}
Several notions of families have been associated with real or complex reflection groups: 
they all form a partition of $\Irr(W)$ and it is a natural question to understand, or compare, 
all these partitions. Fix $c \in \CC$ and let $k=\kappa(c) \in \KC$. Let
$k^\sharp=(k_{\O,j}^\sharp)_{(\O,j) \in \aleph(W)}$ denote the element of $\KC$
defined by $k_{\O,j}^\sharp=k_{\O,-j}$ (the indices $j$ being viewed modulo $e_\O$).
\begin{itemize}
\item[$\bullet$] If $W$ is real (i.e. a Coxeter group) and $c$ has real values 
(since reflections have order $2$ in this case, this is equivalent 
to requiring that $k$ has real values; note also that $k=k^\sharp$ in this case), 
there are two different notions:
\begin{itemize}
\itemth{L} The {\it Lusztig $c$-families} are defined using the notion of {\it $c$-constructible characters} 
of $W$ built on the Lusztig $\ab$-function (see for instance~\cite[\S{15.3.B}]{calogero}). 
They have been explicitly computed by Lusztig in all cases~\cite{lusztig}. Note that they 
are related to the representation theory of finite reductive groups~\cite{lusztig orange}. 

\itemth{KL} The {\it Kazhdan--Lusztig $c$-families} are associated with the partition of $W$ 
into {\it Kazhdan--Lusztig $c$-cells}, built from the 
{\it Kazhdan--Lusztig basis} of the Hecke algebra (see for instance~\cite[\S{8.6}]{calogero}).
\end{itemize}
Lusztig conjectures that both notions coincide~\cite{lusztig} and this conjecture 
holds in many cases (for instance, if $c$ is constant or if $W$ is of type $F_4$). 

\item[$\bullet$] If $W$ is a general complex reflection group, there are also two notions:
\begin{itemize}
\itemth{R} Brou\'e--Kim~\cite{broue kim}
have associated a partition of $\Irr(W)$ into {\it Rouquier $k$-families}
as blocks of the Hecke algebras of $W$ with parameter $k$ over a suitable
ring  (see for instance~\cite[Def.~6.5.1]{calogero} for an extension of this definition to the 
case where $k$ takes complex values: they are called Hecke $k$-families there). 
Chlouveraki~\cite{maria, chlouveraki, chlouveraki B, chlouveraki LNM, chlouveraki D} computed 
Rouquier $k$-families in all cases: as a consequence, she proved that, for $W$ real, 
Rouquier $k$-families coincide with Lusztig $c$-families.
\itemth{CM} The {\it Calogero--Moser $k$-families} defined by Gordon~\cite{gordon} and studied 
in this paper. Note that this last notion of families does not involve the Hecke algebra, 
in contrast to the first three ones. 
\end{itemize}
\end{itemize}

\bigskip
\bigskip

\begin{conjecture}[Martino]\label{conj:martino}
Let $c \in \CC$ and let $k=\kappa(c) \in \KC=\CC$. Then
any Calogero--Moser $c$-family is a union of Rouquier $k^\sharp$-families.
\end{conjecture}

\bigskip
\bigskip

\begin{conjecture}[Gordon--Martino]\label{conj:gordon-martino}
Assume that $W$ is a Coxeter group and assume $c \in \CC$ is real-valued.
Then Calogero--Moser $c$-families coincide with Lusztig $c$-families.
\end{conjecture}

\bigskip

The Martino conjecture is known to be true for the whole infinite series of complex reflection groups $G(m,p,n)$, see \cite{martino, martino 2}. The Gordon--Martino conjecture is known to be true for all Coxeter groups except possibly $H_4$ and $E_6$ -- $E_8$ by \cite{EG, gordon martino, gordon, martino 2, thiel-counter}. 
Therefore, Rouquier families and Calogero--Moser families are two candidates for extending 
the notion of (Kazhdan-)Lusztig families to general complex reflection groups. However, it must 
be pointed out that Martino's Conjecture does not claim that Calogero--Moser families 
coincide with Rouquier families: indeed, there are cases where it is known that both notions 
do not coincide~\cite{thiel-counter}. 

\medskip

We will show that both conjectures are true for all the exceptional complex reflection groups for which we could compute the Calogero--Moser families (see Section \ref{sec:families_summary}). Note that both conjectures involve an {\it infinite} range of parameters. To be able to algorithmically test the conjectures, we need to reduce this to a {\it finite} problem. This is possible due to the {\it semi-continuity} property of both types of families: in both cases there is a finite hyperplane arrangement (the Calogero--Moser hyperplanes for the Calogero--Moser families and the so-called {\it essential hyperplanes} for the Rouquier families \cite{chlouveraki LNM}) such that for all parameters outside these hyperplanes the respective families are always the same and for an arbitrary parameter the families are the smallest unions of the generic families on the hyperplanes containing the parameter. For the Calogero--Moser families this is precisely \eqref{CM_semicont} and for the Rouquier families this property was established in~\cite{calogero}. We thus arrive at the following algorithm:

\bigskip
\refstepcounter{theo}
\noindent{\bfit Algorithm~\thetheo.\label{algo:martino} Testing the Martino conjecture}
\begin{leftbar}
\begin{itemize}
\itemth{1} Test if the generic Rouquier families are unions of the generic Calogero--Moser families.
\itemth{2} For each Calogero--Moser hyperplane defined by a linear polynomial $H$ in the variables $k_{\O,j}$ consider the polynomial $H^\sharp$ obtained by replacing $k_{\O,j}$ by $k_{\O,j}^\sharp = k_{\Omega,-j}$. If $H^\sharp$ is an essential hyperplane, test whether the generic Rouquier families on $H^\sharp$ are unions of the generic Calogero--Moser families on $H$; otherwise test if the generic Rouquier families are unions of the generic Calogero--Moser families on $H$.
\end{itemize}
If all tests are true, then the Martino conjecture is true (for all parameters).
\end{leftbar}

\bigskip

By testing for equality instead of union we can also test the Gordon--Martino conjecture in the same way. We note that the Rouquier families and the essential hyperplanes for the exceptional complex reflection groups have been computed in \cite{chlouveraki LNM} and this data has been added to the database of CHAMP.

\bigskip

\begin{exemple}
We test the Martino conjecture for the group $G_4$ in CHAMP:
\begin{lstlisting}
> W:=ComplexReflectionGroup(4);
> rou := RouquierFamilies(W);
> Keys(rou); //The essential hyperplanes
{
k_{1,2},
k_{1,1} - 2*k_{1,2},
k_{1,1},
2*k_{1,1} - k_{1,2},
1,
k_{1,1} + k_{1,2},
k_{1,1} - k_{1,2}
}
> MartinoConjecture(W);
true
\end{lstlisting}
\end{exemple}

\bigskip

Using this algorithm we have proven:

\bigskip

\begin{theo}\label{theo:gordon-martino} \hfill
  \begin{itemize}
  \itemth{1} The Martino conjecture is true for all exceptional complex reflection groups where we could compute the Calogero--Moser families, i.e. for all   except possibly $G_{16}$ -- $G_{22}$ and $G_{29}$ -- $G_{37}$.

  \itemth{2} The Gordon--Martino conjecture is true for $G_{23} = H_3$ and $G_{28} = F_4$. \qed
\end{itemize}
\end{theo}

\bigskip
\subsection{Cuspidal families}
Recall from Section \ref{sec:symp_sing} that the zero-dimensional symplectic leaves of $\ZC_c$ are also called {\it cuspidal} points. It was shown in \cite{bellamy thiel} that the cuspidal points are contained in $\ZC_c^{\CM^\times}$. Accordingly, we call a Calogero--Moser family {\it cuspidal} if the corresponding $\CM^\times$-fixed point in $\ZC_c$ is cuspidal. The cuspidal families are important because there is a kind of {\it Harish-Chandra theory}, developed by Bellamy~\cite{bellamy cuspidal},
which reduces the study of Calogero--Moser families to the cuspidal ones (for parabolic subgroups of $W$).

For a Calogero--Moser $c$-family $\FC$ we denote by $\Omeb_\FC^c$ the map $\Omeb_\chi^c$ for one (any) $\chi \in \FC$. The corresponding point $\mG_\FC^c$ of $\ZC_c^{\CM^\times}$ is then equal to the kernel of $\Omeb_\FC^c$. If $(z_i)_{i \in I}$ denotes a system of algebra generators of $\Zb$ then $\mG_\FC^c$ is generated by the elements $(z_i-\Omeb_\FC^c(z_i))_{i \in I}$. We thus obtain:

\bigskip

\begin{prop}\label{prop:cuspidal}
$\FC$ is cuspidal if and only if $\Omeb_\FC^c(\{z_i,z_j\}))=0$ for all $i$, $j \in I$.
\end{prop}


Note that, thanks to Lemma~\ref{lem:zplus}, we only need to compute the Poisson
brackets $\{z_i,z_j\}$ for the pairs $(i,j)$ such that $\deg_\ZM(z_i)+\deg_\ZM(z_j)=0$. We thus arrive at the following algorithm:


\bigskip
\refstepcounter{theo}
\noindent{\bfit Algorithm~\thetheo.\label{algo:cuspidal}  Cuspidal families}
\begin{leftbar}
Let $c \in \CC$.
\begin{itemize}
\itemth{1} Compute a system of generators  $(z_i)_{i \in I}$ of $\Zb$ via  Algorithm~\ref{algo:generateurs}.
\itemth{2} Compute the Calogero--Moser $c$-families via  Algorithm~\ref{algo:familles}.
\itemth{3} Compute the Poisson brackets $\{z_i,z_j\}_{i,j \in I}$ for all $(i,j)$ such that $\deg_\ZM(z_i)+\deg_\ZM(z_j)=0$ via Algorithm~\ref{algo:brackets}.
\end{itemize}
Then the cuspidal Calogero--Moser $c$-families are precisely those families $\FC$ for which $\Omeb_\FC^c(\{z_i,z_j\}))=0$ for all $(i,j)$ such that $\deg_\ZM(z_i)+\deg_\ZM(z_j)=0$.
\end{leftbar}

\bigskip

\begin{exemple}
We compute the cuspidal Calogero--Moser families for the dihedral group of order 8 in CHAMP and verify that in the equal parameter case there is a unique cuspidal family (see also Section \ref{sec:dihedral_groups}):
\begin{lstlisting}
> W:=ComplexReflectionGroup(4,4,2);
> H:=RationalCherednikAlgebra(W,0);
> cusp := CuspidalCalogeroMoserFamilies(H);
> R := Universe(cusp);
> cusp[R.1-R.2];
{
{ 1, 2, 5 }
}
\end{lstlisting}
\end{exemple}

\bigskip

\section{Cellular characters}

\medskip
\def\gaudin{{\mathrm{Gau}}}

We fix $c \in \CC$. Calogero--Moser $c$-cellular characters have been defined
in~\cite[Def.~1.8.4]{calogero} but we will use here the equivalent definition
given by~\cite[Theo.~8.3.2]{calogero} in terms of Gaudin algebras.

\bigskip

\subsection{Definition}
We recall here the definition of cellular characters which involves
the {\it Gaudin algebra}~\cite[\S{13.2}]{calogero}.
First, recall that $\kb[\CC \times V^\reg][W]$
denotes the group algebra of $W$ over the algebra $\kb[\CC \times V^\reg]$, and not the semi-direct
product $\kb[\CC \times V^\reg] \rtimes W$. For $y \in V$, let
$$\DC_y=\sum_{s \in \REF(W)}
\e(s)C_s \frac{\langle y,\a_s\rangle}{\a_s} s\indexnot{D}{\DC_y}\,\,\in \kb[\CC \times V^\reg][W].$$
Then $\gaudin(W)$ is the sub-$\kb[\CC \times V^\reg]$-algebra of
$\kb[\CC \times V^\reg][W]$ generated by
the $\DC_y$'s ($y \in V$): it will be called the {\it generic Gaudin algebra}
associated with $W$. Recall that it is commutative~\cite[\S{8.3.B}]{calogero}.

Now, if $y \in V$, let $\DC_y^c \in \kb[V^\reg][W]$ denote the specialization at $c$ of $\DC_y$
and let $\gaudin_c(W)$ denote the specialization of $\gaudin(W)$ at $c \in \CC$.
Since $\kb(V)^W \otimes_{\kb[V^\reg]^W} \kb[V^\reg] \simeq \kb(V)$, we have
$$\CM(V)^W \gaudin(W)=\kb(V)\gaudin_c(W)$$
and $\kb(V)\gaudin_c(W)$ is the $\kb(V)$-sub-algebra of $\kb(V)[W]$ generated by the
$\DC_y^c$'s.

Let $L$ be a simple $\kb(V)\gaudin_c(W)$-module and set
$$\g_L^\gaudin = \sum_{\chi \in \Irr(W)} [\Res_{\kb(V)\gaudin_c(W)}^{\kb(V)[W]} \kb(V) \otimes V_\chi : L] ~\chi.$$
Then $\g_L^\gaudin$ is a Calogero--Moser $c$-cellular character of $W$,
and they can all be obtained in this way~\cite[Theo.~13.2.10]{calogero}.
Note moreover that $W$ acts on $\gaudin_c(W)$, so it
acts on $\Irr(\kb(V)\gaudin_c(W))$ and
\equat\label{eq:conjugaison-cellular}
\g_{\lexp{w}{L}}^\gaudin=\g_L^\gaudin
\endequat
for all $L \in \Irr(\kb(V)\gaudin_c(W))$ and all $w \in W$.

\bigskip
\def\calopetit{{\SSS{\mathrm{CM}}}}

\bigskip

\def\kM{\Bbbk}

\subsection{Algorithm}
Following Appendix~\ref{appendix:cellulaire} (see Remark~\ref{comment:algo-cellular}), the following
algorithm allows us to compute the Calogero--Moser $c$-cellular characters. For this, we fix a number field
$\kM$ and a $\kM$-form $V_\kM$ of $V$ which is stabilized by $W$ (this is done in order
to perform exact computations).

\bigskip

\refstepcounter{theo}
\noindent{\bfit Algorithm~\thetheo.\label{algo:cellulaires}  Cellular characters}
\begin{leftbar}
Let $c \in \CC$ having values in $\kM$.
The set of Calogero--Moser $c$-cellular characters is computed as follows:
\begin{itemize}
\itemth{1} Let $\DC \in \End_{\kM[V_\kM \times V_\kM^\reg]}(\kM[V_\kM \times V_\kM^\reg][W])$
be the element such that the action by left translation of $\DC_y$ on $\kM[V_\kM^\reg][W]$
is the specialization at $y$ of $\DC$ (for all $y \in V$).
\itemth{2} Let $\Pi_\DC$ denote the quotient of the characteristic polynomial $\Pit_\DC$
of $\DC$ by the greatest common divisor of $\Pit_\DC$ and $\Pit_\DC'$. Let $\D_\DC$
denote its discriminant (it is an element of $\kM[V_\kM \times V_\kM^\reg]$). Find (randomly or
algorithmically) $(y_\kM,v) \in V_\kM \times V_\kM^\reg$ such that $\D_\DC(y_\kM,v) \neq 0$.
\itemth{3} Let $\DC_y^v$ denote the specialization of $\DC$ at $(y,v)$ and let
$\Pi_{y_\kM}^v$ denote the specialization of $\Pi$ at $(y_\kM,v)$. Compute the
factorization $\Pi_{y_\kM}^v=\Pi_1\cdots\Pi_r$ into a product of irreducible
polynomials in $\kM[\tb]$.
\itemth{4} For $1 \le i \le r$, compute the character $\G_i$ of the representation of
$W$ acting on $\Ker(\Pi_i(\DC_{y_\kM}^v)^{n_i})$ (where $n_i$ is the valuation of $\Pi_i$
in the characteristic polynomial of $\DC_{y_\kM}^v$).
\itemth{5} By Remark~\ref{comment:algo-cellular},
$\{\G_i/\deg(\Pi_i)~|~1 \le i \le r\}$ is the set
of $c$-cellular characters of $W$.
\end{itemize}
\end{leftbar}

\bigskip

The fact that the set of Calogero--Moser $c$-cellular characters
does not depend on the choice of the number field $\kM$
follows from Remark~\ref{rem:extension-field} and from the fact that the group algebra $\kM W$
is split~\cite{benard, bessis}.

\bigskip

\begin{rema}\label{comment:reps}
The matrix which gives $\DC$ has size $|W|$, so the computations involved in the
previous algorithm can be heavy even for reasonably small $W$'s (type $H_3$, $F_4$, $H_4$\dots).
If one has all the explicit representations of $W$, then the computation of $\Pi$
might be performed for each representation, instead of the regular representation.\finl
\end{rema}

\bigskip

\subsection{Spetses}\label{sub:spetses}
For $W$ a Coxeter group, it is conjectured that Calogero--Moser cellular characters
coincide with Lusztig constructible characters~\cite[Conj.~L]{calogero}. This has
been checked for dihedral groups~\cite{bonnafe diedral} and for the symmetric group~\cite{cells-typeA},
as well when the Calogero--Moser space is smooth~\cite[Theo.~14.4.1]{calogero}.

\medskip

For $W$ a {\it spetsial} complex reflection group and for the {\it spetsial} parameter $c_{\srm\prm}$ 
(i.e. the parameter such that $\kappa(c_{\srm\prm})_{\O,0}=1$ and $\kappa(c_{\srm\prm})_{\O,j}=0$ 
for all $\O \in \AC/W$ 
and $1 \le j \le e_\O-1$),
the philosophy of {\it spetses}~\cite{BMM} allows one to propose good candidates for
a notion of constructible characters~\cite{malle rouquier} (let us call
them Malle--Rouquier $c_{\srm\prm}$-constructible characters). The construction
is completely different from the construction of cellular characters, but it
is conjectured that both notions coincide. As explained above, computations
can become pretty heavy when the size of the group increases, but we were able to check 
the following result:

\bigskip

\begin{theo}
If $W$ a spetsial primitive complex reflection group of rank $\le 2$
(i.e. if $W=G_4$, $G_6$, $G_8$ or $G_{14}$), then Calogero--Moser 
$c_{\srm\prm}$-cellular characters coincide with Malle--Rouquier 
$c_{\srm\prm}$-constructible characters.
\end{theo}

\bigskip

Even though it seems like a weak result, it is an extension of the
surprising links between the Poisson geometry of Calogero--Moser spaces
and unipotent representations of finite reductive groups explained by the first author
in~\cite{cm-unip}.

\bigskip

\setcounter{section}{0}
\renewcommand\thesection{\Roman{section}}
\def\sectionname{Appendix}

\section{A general theory of cellular characters and how to compute them}\label{appendix:cellulaire}

If $A$ is a commutative ring, $M$ is free $A$-module of finite rank and
$\ph \in \End_A(M)$, we denote by $\carac_\ph \in A[\tb]$ the characteristic
polynomial of $\ph$. If $K$ is a field and
$f \in K[\tb]$, we denote by $f^\rad$ the greatest common divisor of $f$ and $f'$
(which is normalized to be monic) and we set $f^\semi=f/f^\rad \in K[\tb]$.
If moreover $A=K$, we denote by $\mini_\ph \in K[\tb]$
the minimal polynomial of $\ph$ (which is also normalized to be monic).
Note that $\carac_\ph^\semi = \mini_\ph^\semi$.

\medskip

Let $\kM$ be a field of characteristic $0$.

\bigskip

\boitegrise{{\bf Hypothesis and notation.}
{\it We fix in this appendix a finite dimensional $\kM$-vector space $E$,
a {\bfit split} subalgebra $A$ of $\End_\kM(E)$ and an integral and integrally closed
$\kM$-algebra $P$. Its field of fraction is denoted by $K$.\\
\hphantom{aa}If $R$ is a commutative $\kM$-algebra, we set $RE=R \otimes E$ and
$RA=R \otimes A$. We also fix a positive integer $n$ and $n$ commuting elements
$D_1$,\dots, $D_n$ of $\End_P(PE)$ (for $1 \le i \le n$) which also
commute with the action of $PA$.\\
\hphantom{aa}Finally, we set $D=(D_1,\dots,D_n)$ and we denote by $P[D]$ the commutative
algebra generated by $D_1$,\dots, $D_n$. If $\pG$ is a prime ideal
of $P$, we denote by $D_i(\pG) \in \End_{P/\pG}((P/\pG)E)$ the reduction of $D_i$ modulo $\pG$
and we set $D(\pG)=(D_1(\pG),\dots,D_n(\pG))$.}}{0.75\textwidth}

\bigskip

We aim to define a notion of {\it cellular characters} of $A$ and study its behavior under
specialization or field extension. The results presented here are certainly neither new nor
original: they are just adapted to our situation, and stated here with proof for the sake
of completeness.

\bigskip

\subsection{Cellular characters}
The $K$-algebra $KA$ is split, so there is a bijection~\cite[Propositions~3.56~and~7.7]{curtis-reiner}
$$\fonctio{\Irr(K[D]) \times \Irr(KA)}{\Irr(K[D] \otimes_{K} KA)}{(L_1,L_2)}{
L_1 \otimes L_2.}$$
This induces an isomorphism of $\ZM$-modules
$$\groth(K[D]) \otimes_\ZM \groth(KA) \longiso
\groth(K[D] \otimes_K KA),$$
where, for a finite dimensional algebra $B$, 
$\groth(B)$ is the Grothendieck group of the category of $B$-modules. 
But the $K$-vector space $KE$ inherits an action of
$K[D] \otimes_K KA$. So its class $\isomorphisme{KE}$
in $\groth(K[D] \otimes_K KA)$ can be written
\equat\label{eq:classe-E}
\isomorphisme{KE} = \sum_{S \in \Irr(K[D])} \isomorphisme{S} \otimes \g_S^P,
\endequat
where $\g_S^p \in \groth(KA) \simeq \groth(A)$ (because $A$ is split).

\bigskip

\begin{defi}\label{defi:cellulaire-D}
The element $\g_S^P$ of $\groth(A)$ is called a {\bfit $D$-cellular character}.
\end{defi}

\bigskip

It might be the case that the map
$$\fonction{\g^P}{\Irr(K[D])}{\groth(A)}{S}{\g_S^P}$$
is not injective. But we will only be interested in the image of $\g^P$,
and its behavior under specialization or field extension.
This last situation is the simplest:

\bigskip

\begin{prop}\label{prop:field-extension}
Let $\kM'$ be an extension field of $\kM$ and let $P'=\kM' \otimes P$.
We assume that $P'$ is integral and integrally closed and we set $K'=\kM' \otimes K$
(it is then the fraction field of $P'$). Then the images
of the maps $\g^P$ and $\g^{P'}$ are equal.
\end{prop}

\bigskip

\begin{proof}
Note that $K[D]$ is commutative. So, if $S$ is a simple
$K[D]$-module, then $\kM' \otimes S$ is a multiplicity
free semisimple $K'[D]$-module~\cite[Propositions~3.56~and~7.7]{curtis-reiner}.
Moreover, $\kM' \otimes S_1$ and $\kM' \otimes S_2$ have a common irreducible constituent if and only if
$S_1 \simeq S_2$ (see~\cite[Theorem~7.9]{curtis-reiner}). The proposition follows.
\end{proof}

\bigskip

\subsection{The case ${\boldsymbol{n=1}}$}
We will first study here the case where $n=1$. As it will be explained in the next
subsection, the general case can in fact be reduced to this case.

\bigskip

\boitegrise{\noindent{\bf Hypothesis and notation.}
{\it We assume in this subsection, and only in this subsection,
that $n=1$ and we write $D=D_1$.}}{0.75\textwidth}

\bigskip

The characteristic polynomial $\carac_D$ of $D$ has coefficients in $P$. Using Euclid's algorithm
in $K[\tb]$, one can compute the polynomial $\carac_D^\rad$ (recall that it is
the greatest common divisor of $\carac_D$ and its derivative $\carac_D'$)
and so the polynomial $\carac_D^\semi=\carac_D/\carac_D^\rad$ is computable. We then write
$$\carac_D^\semi=\Pi_1 \cdots \Pi_r,\qquad \carac_D=\Pi_1^{n_1}\cdots \Pi_r^{n_r}$$
where $\Pi_i$ belongs to $K[\tb]$ and is irreducible, $n_i \ge 1$ and $\Pi_i \neq \Pi_j$
if $i \neq j$. By~\cite[Proposition~B.4.1]{calogero},
\equat\label{eq:pi-kx}
\carac_D^\semi,\Pi_1,\dots,\Pi_r \in P[\tb].
\endequat
If $1 \le i \le r$, we set
$$\LC_i=K[\tb]/\langle \Pi_i \rangle,$$
and we view it as a $P[D]$-module, where $D$ acts by multiplication by
the indeterminate $\tb$. Then
\equat\label{eq:simple-kd}
\Irr(K[D]) = \{\LC_1,\dots,\LC_r\}\qquad \text{and} \qquad
\text{$\LC_i \not\simeq \LC_j$ if $i \neq j$.}
\endequat
Then
\equat\label{eq:lemme-noyau}
KE = \bigoplus_{i=1}^r \Ker(\Pi_i(D)^{n_i})
\endequat
is a decomposition of $KE$ into $K[D] \otimes_K KA$-modules. As a $K[D]$-module,
the chief factors of $\Ker(\Pi_i(D)^{n_i})$ are all isomorphic to $\LC_i$. Hence,
if we denote by $\isomorphisme{\Ker(\Pi_i(D)^{n_i})}_{KA}$ the image of
$\Ker(\Pi_i(D)^{n_i})$ in $\groth(KA)$, then this image is a multiple of
$\dim(\LC_i)=\deg(\Pi_i)$ and
\equat\label{eq:cellulaires}
\isomorphisme{KE}=\sum_{i=1}^r \isomorphisme{\LC_i} \otimes_\ZM \Bigl(\frac{1}{\deg(\Pi_i)}
\isomorphisme{\Ker(\Pi_i(D)^{n_i})}_{KA}\Bigr).
\endequat
Therefore:

\bigskip

\begin{prop}\label{prop:cellulaires-generaux}
With the above notation, the set of $D$-cellular characters of $A$ is equal to
$$\Bigl\{\frac{1}{\deg(\Pi_i)}
\isomorphisme{\Ker(\Pi_i(D)^{n_i})}_{KA}~\Bigl|~1 \le i \le r\Bigr\}.$$
\end{prop}

\bigskip

%
%
%
%
%

We denote by $\D$ the discriminant of the polynomial $\carac_D$: then $\D \in P$ and, if $\pG$
is a prime ideal of $P$, we denote by $\D(\pG)$ the image of $\D$ in $P/\pG$. Note that
\equat\label{eq:disc-non-zero}
\D \neq 0
\endequat
because the factorization of $\carac_D$ into irreducible polynomials is multiplicity-free
(and since $\kM$ has characteristic $0$).

\bigskip

\begin{prop}\label{prop:specialization}
Let $\pG$ be a prime ideal of $P$ such that $\D(\pG) \neq 0$ and $P/\pG$ is
integrally closed. Then the set of $D$-cellular characters coincides
with the set of $D(\pG)$-cellular characters.
\end{prop}

\bigskip

\begin{proof}
Assume that $\D(\pG) \neq 0$. Let us write
$$\Pi_i(\pG)=\prod_{j=1}^{d_i} \pi_{i,j}^{e_{i,j}},$$
where $\pi_{i,j} \in \kM_P(\pG)[\tb]$ is irreducible, $e_{i,j} \ge 1$ and $\pi_{i,j} \neq \pi_{i,j'}$
if $1 \le j < j' \le d_i$. Then
$$\carac_D(\pG) = \prod_{i=1}^r \prod_{j=1}^{d_i} \pi_{i,j}^{e_{i,j}}.$$
Since $\D(\pG) \neq 0$, this means that $e_{i,j} = 1$ for all $i$, $j$ and that
$\pi_{i,j}=\pi_{i',j'}$ if and only if $(i,j)=(i',j')$.

Now, let $\EC=\{(i,j)~|~1\le i \le r$ and $1 \le j \le d_i\}$ and, if $(i,j) \in \EC$, let
$$\LC_{i,j}^\pG=\kM_P(\pG)[\tb]/\langle \pi_{i,j} \rangle,$$
viewed as a $\kM_P(\pG)[D(\pG)]$-module, where $D(\pG)$ acts by multiplication by $\tb$.
Then, by~(\ref{eq:simple-kd}), the map
$$\fonctio{\EC}{\Irr(\kM_P(\pG)[D(\pG)])}{(i,j)}{\LC_{i,j}^\pG}$$
is bijective. We only need to show that
$$\g_{\LC_i}^P=\g_{\LC_{i,j}^\pG}^{P/\pG}.\leqno{(*)}$$
Since $\kM$ has characteristic $0$, an equality between elements of the Grothendieck group
is equivalent to an equality between the corresponding (virtual) characters.
If we denote by $\chi_i$ (resp. $\chi_{i,j}^\pG$) the
character of $\LC_i$ (resp. $\LC_{i,j}^\pG$), then
$$\chi_i \equiv \sum_{j=1}^{d_i} \chi_{i,j} \mod \pG.$$
Since $\chi_{i,j}=\chi_{i',j'}$ if and only if $(i,j)=(i',j')$,
the result follows from~(\ref{eq:classe-E}) by taking associated characters
and reduction modulo $\pG$.
\end{proof}

\bigskip

%
%
%
%

\bigskip

\subsection{Back to the general case}
We come back to the initial set-up of this appendix.

\bigskip

\boitegrise{{\bf Hypothesis and notation.} {\it Until the end of this appendix,
we no longer assume that $n=1$.
}}{0.75\textwidth}

\bigskip

The main result of this subsection (see Proposition~\ref{prop:cellular-n}) shows that
the computation of the $D$-cellular characters can be reduced to the case where $n=1$,
by extending the ring of definition $P$. Let $X_1$,\dots, $X_n$ be indeterminates over $P$
and let $P[\Xb]=P[X_1,\dots,X_n]$. We set
$$\DC = X_1 D_1 + \cdots + X_n D_n \in \End_{P[\Xb]}(P[\Xb]E).$$
Since $P[\Xb]$ is integrally closed, we are in the same set-up as in the previous
section (the case $n=1$), with $P$ replaced by $P[\Xb]$ and $D$ replaced by $\DC$.
So we can define similarly $\DC$-cellular characters of $A$.

\bigskip

\begin{prop}\label{prop:cellular-n}
The set of $D$-cellular characters of $A$ coincides with the set
of $\DC$-cellular characters of $A$.
\end{prop}

\bigskip

\begin{proof}
Let $\D_\Xb$ denote the discriminant of $\carac_\DC^\semi$ (it is an element of $P[\Xb]$).
Since $\D_\Xb \neq 0$, there exists $x=(x_1,\dots,x_n) \in \kM^n$ such that
$$\D_\Xb(x_1,\dots,x_n) \neq 0.$$
Let $\DC(x)=x_1 D_1 + \cdots + x_n D_n \in \End_P(PE)$. It follows from
Proposition~\ref{prop:specialization} that the set of $\DC$-cellular characters
coincide with the set of $\DC(x)$-cellular characters (one must replace, in Proposition~\ref{prop:specialization},
the ring $P$ by $P[\Xb]$ and take $\pG=\langle X_1-x_1,\dots,X_n-x_n\rangle$). So it remains to prove that
the set of $\DC(x)$-cellular characters coincides with the set of $D$-cellular characters.

Since $K[\DC(x)] \subset K[D]$, it remains to prove that $K[D]=K[\DC(x)] + \Rad(K[D])$.
By~\cite[Corollary~7.8~and~Theorem~7.9]{curtis-reiner}, we may assume that $K$
is algebraically closed. In this case, there exists a basis of $KE$
such that $D_1$,\dots, $D_n$ are represented by commuting upper triangular matrices.
Let $d=\dim_\kM(E)$ and let $(\l_j^{(i)})_{1 \le j \le d}$ denote the
(ordered) sequence of diagonal coefficients of $D_i$ and let
$$\fonction{\lamb}{\{1,2,\dots,d\}}{\kM^n}{j}{(\l_j^{(1)},\dots,\l_j^{(n)}).}$$
Let $\EC$ denote the image of the map $\lamb$.
Then
$$\carac_{\DC}^\semi=\prod_{(\mu_1,\dots,\mu_n) \in \EC} (\tb-\sum_{i=1}^n X_i \mu_i)$$
and the condition $\D_\Xb(x) \neq 0$ implies that
$$\carac_{\DC(x)}^\semi=\prod_{(\mu_1,\dots,\mu_n) \in \EC} (\tb-\sum_{i=1}^n x_i \mu_i)$$
and that
$$\text{$\sum_{i=1}^n x_i \mu_i=\sum_{i=1}^n x_i \mu_i'$ if and only if
$(\mu_1,\dots,\mu_n)=(\mu_1',\dots, \mu_n')$.}\leqno{(*)}$$

Now, let $M$ be an element of $K[D]$, and let $m_j$ denote its $j$-th diagonal coefficient.
If $j$, $j'$ are such that $\lamb(j)=\lamb(j')$, then necessarily $m_j=m_{j'}$.
It then follows from $(*)$ that there exists a polynomial $f \in K[\tb]$ such that
$$f(\sum_{i =1}^n x_i \l_j^{(i)})=m_j$$
for all $j \in \{1,2,\dots,d\}$. Therefore, $M-f(\DC(x))$ belongs to $K[D]$
and has all its diagonal coefficients equal to $0$. This forces
$M-f(\DC(x)) \in \Rad(K[D])$, as desired.
\end{proof}

\bigskip

\begin{rema}\label{comment:algo-cellular}
Proposition~\ref{prop:cellular-n} together with Proposition~\ref{prop:cellulaires-generaux}
provides a way to reduce the problem of computing $D$-cellular characters to the case $n=1$.
In this former case, using Proposition~\ref{prop:specialization}, one can reduce
the problem to the case where $D=D_1$ has coefficients in a finite extension of $\kM$
(if $P$ is finitely generated over $\kM$).\finl
\end{rema}

\bigskip

\begin{rema}\label{rem:extension-field}
In this remark, and only in this remark, the $D$-cellular characters
will be called the $(\kM,D)$-cellular characters.
Let $\kM'$ be a field extension of $\kM$ such that $P'=\kM' \otimes_\kM P$ is
still an integrally closed domain. Let $K'$ be the fraction field of $P'$.
Since the $\kM'$-algebra $\kM' A$ is still split, we can
define the $(\kM',D)$-cellular characters $\g_{S'}^{P'}$ attached
to a simple $K'[D]$-module $S'$. In fact, through the isomorphism
$\groth(A) \longiso \groth(\kM' A)$ induced by scalar extension,
$$\text{\it the $(\kM,D)$-cellular characters and the $(\kM',D)$-cellular characters coincide.}$$
Indeed, since we are in characteristic $0$ and since $K[D]$ is commutative, any simple
$K[D]$-module $S$ satisfies
$$K' \otimes_K S \simeq \bigoplus_{j=1}^r S_j',$$
where $S_1'$,\dots, $S_r'$ are two-by-two non-isomorphic simple $K'[D]$-modules. Moreover,
if $S_1$ and $S_2$ are two non-isomorphic simple $K[D]$-modules, then
$\Hom_{K'[D]}(K' \otimes_K S_1,K' \otimes_K S_2)=K' \otimes_K \Hom_{K[D]}(S_1,S_2)=0$. Both facts
imply that
$$\g_S^P=\g_{S_1'}^{P'}=\cdots=\g_{S_r'}^{P'},$$
as desired.\finl
\end{rema}

\bigskip


\begin{thebibliography}{AAAA}

\bibitem[Bea]{beauville} {\sc A. Beauville},
{\it Symplectic singularities},
Invent. Math. {\bf 139} (2000), 541-549.

\bibitem[BCHM]{BCHM} {\sc C. Birkar, P. Cascini, C.D. Hacon, J. McKernan},
{\it Existence of minimal models for varieties of log general type},
J. Amer. Math. Soc. {\bf 23} (2010), no. 2, 405–468.

%
%
%
%
\bibitem[Bel1]{bellamy cuspidal} {\sc G. Bellamy},
{\it Cuspidal representations of rational Cherednik algebras at $t=0$},
Math. Z. {\bf 269} (2011), 609-627.

\bibitem[Bel2]{bellamy counting} {\sc G. Bellamy},
{\it Counting resolutions of symplectic quotient singularities},
Compos. Math. {\bf 152} (2016), 99-114.

\bibitem[BBFJLS]{BBFJLS} {\sc G. Bellamy, C. Bonnaf\'e, B. Fu, D. Juteau, M. Levy \& E. Sommers},
{\it A new family of isolated symplectic singularities
with trivial local fundamental group},
Proc. London Math. Soc. {\bf 126} (2023), 1496-1521.

%

\bibitem[BeMaSc]{be-sc} {\sc G. Bellamy, R. Maksimau \& T. Schedler}, in preparation.

\bibitem[BST2]{BST-Hyperplanes} {\sc G. Bellamy, T. Schedler \& U. Thiel},
{\it Hyperplane arrangements associated with symplectic quotient singularities},
Phenomenological approach to algebraic geometry, 25–45, Banach Center Publ., 116, Polish Acad. Sci. Inst. Math., Warsaw, 2018.

\bibitem[BST1]{BST-Towards} {\sc G. Bellamy, J. Schmitt \& U. Thiel},
{\it Towards the classification of symplectic linear quotient singularities admitting a symplectic resolution}, to Math. Z. {\bf 300} (2022), 661-681.

\bibitem[BelThi]{bellamy thiel} {\sc G. Bellamy \& U. Thiel},
{\it Cuspidal Calogero--Moser and Lusztig families for Coxeter groups},
J. Algebra {\bf 462} (2016), 197-252.


\bibitem[Bena]{benard} {\sc M. Benard},
{\it Schur indices and splitting fields of the unitary reflection groups},
{\it J. Algebra} {\bf 38} (1976), 318-342.

\bibitem[Bens]{benson} {\sc D. J. Benson},
{\it Polynomial invariants of finite groups},
London Math. Soc. Lecture Note Series {\bf 190}, Cambridge University Press, Cambridge, 1993.

\bibitem[Bes]{bessis} {\sc D. Bessis},
{\it Sur le corps de d\'efinition d'un groupe de r\'eflexions complexe},
Comm. Algebra {\bf 25} (1997), 2703-2716.
%

%

\bibitem[Bon1]{bonnafe diedral} {\sc C. Bonnaf\'e},
{\it On the Calogero--Moser space associated with dihedral groups},
Ann. Math. Blaise Pascal {\bf 25} (2018), 265-298.

\bibitem[Bon2]{bonnafe diedral 2} {\sc C. Bonnaf\'e},
{\it On the Calogero--Moser space associated with dihedral groups II,
The equal parameter case}, preprint (2021), {\tt arXiv:2112.12401}. To appear
in Ann. Math. Blaise Pascal.

\bibitem[Bon3]{auto} {\sc C. Bonnaf\'e},
{\it Automorphisms and symplectic leaves of Calogero--Moser spaces},
J. Australian Math. Soc. {\bf 115} (2023), 26-57.

\bibitem[Bon4]{regular} {\sc C. Bonnaf\'e},
{\it Regular automorphisms and Calogero--Moser families},
preprint (2021), {\tt arxiv:2112.13685}.

\bibitem[Bon5]{cm-unip} {\sc C. Bonnaf\'e},
{\it Calogero--Moser spaces vs unipotent representations},
preprint (2021), {\tt arXiv:2112.13684}. To appear in Pure and App. Math. Quart.


\bibitem[BoMa]{bonnafe maksimau} {\sc C. Bonnaf\'e \& R. Maksimau},
{\it Fixed points in smooth Calogero--Moser spaces},
Ann. Inst. Fourier {\bf 71} (2021), 643-678.



%
%
%
%
%

%
%
%
%
%
%

\bibitem[BoRo1]{calogero-first} {\sc C. Bonnaf\'e \& R. Rouquier},
{\it Calogero-{M}oser versus {K}azhdan-{L}usztig cells},
Pacific J. Math. {\bf 261} (2013), 45-51.

\bibitem[BoRo2]{calogero} {\sc C. Bonnaf\'e \& R. Rouquier},
{\it Cherednik algebras and Calogero--Moser cells},
{\tt arXiv:1708.09764}.

\bibitem[Bou]{bourbaki} {\sc N. Bourbaki},
{\it Alg\`ebre commutative, chapitres 5, 6, 7}.

%
\bibitem[Bro]{broue} {\sc M. Brou\'e},
{\it Introduction to complex reflection groups and their braid groups},
Lecture Notes in Mathematics {\bf 1988}, 2010, Springer.



\bibitem[BrGoWh]{cells-typeA} {\sc A. Brochier, I. Gordon \& N. White},
{\it Gaudin Algebras, RSK and Calogero--Moser Cells in Type A},
Proc. Lond. Math. Soc. {\bf 126} (2023), 1467-1495.

\bibitem[BrKi]{broue kim} {\sc M. Brou\'e \& S. Kim},
{\it Familles de caract\`eres des alg\`ebres de Hecke cyclotomiques},
Adv. Math. {\bf 172} (2002), 53-136.

\bibitem[BrMaMi1]{BMM} {\sc M. Brou\'e, G. Malle \& J. Michel},
{\it Towards Spetses I}, Transf. Groups {\bf 4} (1999) 157-218.

%
%

\bibitem[BrGo]{BrGo} {\sc K. Brown \& I. Gordon},
{\it Poisson orders, symplectic reflection algebras and representation theory},
J. reine angew. Math. {\bf 559} (2003), 193—216.

%
%
%

\bibitem[Cal]{CalogeroOriginal}{\sc F. Calogero},
{\it Solution of the one-dimensional $n$-body problems with quadratic and/or inversely quadratic pair potentials}, J. Math. Phys. {\bf 12} (1971), 419--436

\bibitem[Chl1]{maria} {\sc M. Chlouveraki},
{\it Sur les alg\`ebres de Hecke cyclotomiques des groupes de r\'eflexions complexes},
th\`ese, Paris 7 (2007).

\bibitem[Chl2]{chlouveraki} {\sc M. Chlouveraki},
{\it Rouquier blocks of the cyclotomic Hecke algebras},
C. R. Math. Acad. Sci. Paris {\bf 344} (2007), 615-620.

\bibitem[Chl3]{chlouveraki B} {\sc M. Chlouveraki},
{\it Rouquier blocks of the cyclotomic Ariki-Koike algebras},
Algebra Number Theory {\bf 2} (2008), 689-720.

\bibitem[Chl4]{chlouveraki LNM} {\sc M. Chlouveraki},
{\it Blocks and families for cyclotomic Hecke algebras},
Lecture Notes in Mathematics {\bf 1981}, 2009, Springer.

\bibitem[Chl5]{chlouveraki D} {\sc M. Chlouveraki},
{\it Rouquier blocks of the cyclotomic Hecke algebras of $G(de,e,r)$},
Nagoya Math. J. {\bf 197} (2010), 175-212.

\bibitem[CuRe]{curtis-reiner} {\sc C. W. Curtis \& I. Reiner},
{\it Methods of representation theory, Vol. I, With applications to finite groups and orders}.
Pure and Applied Mathematics, A Wiley-Interscience Publication, John Wiley \& Sons,
Inc., New York, 1981.

\bibitem[EtGi]{EG} {\sc P. Etingof \& V. Ginzburg}, Symplectic reflection algebras,
Calogero--Moser space, and deformed Harish-Chandra homomorphism,
{\it Invent. Math.} {\bf 147} (2002), 243-348.

%
%
%
%
%
%
%
\bibitem[GGOR]{ggor} {\sc V. Ginzburg, N. Guay, E. Opdam \& R. Rouquier},
On the category $\OC$ for rational Cherednik algebras, {\it Invent. Math.} {\bf 154}
(2003), 617-651.

\bibitem[GiKa]{GK} {\sc V. Ginzburg \& D. Kaledin},
Poisson deformations of symplectic quotient singularities,
{\it Adv. in Math.} {\bf 186} (2004), 1-57.

\bibitem[Gor]{gordon} {\sc I. Gordon},
{\it Baby Verma modules for rational Cherednik algebras},
Bull. London Math. Soc. {\bf 35} (2003), 321-336.

\bibitem[Gor]{gordon icra} {\sc I. Gordon},
{\it Symplectic reflection algebras},
in {\it Trends in representation theory of algebras and related topics}, 285-347,
EMS Ser. Congr. Rep., Eur. Math. Soc., Z\"urich, 2008.



\bibitem[GoMa]{gordon martino} {\sc I. G. Gordon \& M. Martino},
{\it Calogero--Moser space, restricted rational Cherednik algebras and two-sided cells},
Math. Res. Lett. {\bf 16} (2009), 255-262.

%

\bibitem[KaLu]{KL} {\sc D. Kazhdan \& G. Lusztig},
{\it Representations of Coxeter groups and Hecke algebras},
Invent. Math. {\bf 53} (1979), 165-184.

\bibitem[Kin]{king}{ \sc S. King},
{\it Minimal generating sets of non-modular invariant rings of finite groups},
J. Symbolic Comput. {\bf 48} (2013), 101--109.


\bibitem[Kol]{kollar}{\sc J. Kollár},
{\it Singularities of the minimal model program},
Cambridge Tracts in Mathematics, 200. Cambridge University Press, Cambridge, 2013.


\bibitem[LeSO]{lehn-sorger} {\sc M. Lehn \& C. Sorger},
{\it A symplectic resolution for the binary tetrahedral group},
Geometric methods in representation theory. II, 429–435,
{\it Sémin. Congr.}, 24-II, Soc. Math. France, Paris, 2012.

\bibitem[Lus]{Lu2} {\sc G.~Lusztig},
	{\it Left cells in Weyl groups},
	in {\it Lie group representations, I}, 99-111,
	Lecture Notes in Math. {\bf 1024}, Springer, Berlin, 1983.

\bibitem[Lus2]{lusztig orange} {\sc G. Lusztig}, {\it Characters of reductive groups over finite fields},
 Ann. Math. Studies {\bf 107}, Princeton UP (1984).

\bibitem[Lus3]{lusztig}
{\sc G.~Lusztig},
\emph{Hecke algebras with unequal parameters}, CRM Monograph Series
{\bf 18}, American Mathematical Society, Providence, RI (2003).

\bibitem[Magma]{magma} {\sc W. Bosma, J. Cannon \& C. Playoust},
{\it The Magma algebra system. I. The user language},
J. Symbolic Comput. {\bf 24} (1997), 235-265.


\bibitem[Mal]{malle icm} {\sc G. Malle},
{\it Spetses}, Proceedings of the International Congress of Mathematicians, 
Vol. II (Berlin, 1998). Doc. Math. 1998, Extra Vol. II, 87-96.

%
\bibitem[MalRou]{malle rouquier} {\sc G. Malle \& R. Rouquier},
Familles de caract\`eres de groupes de r\'eflexions complexes,
{\it Representation Theory} {\bf 7} (2003), 610-640.
%

\bibitem[Mar1]{martino} {\sc M. Martino},
{\it The Calogero--Moser partition and Rouquier families for complex reflection groups},
J. of Algebra {\bf 323} (2010), 193-205.

\bibitem[Mar2]{martino 2} {\sc M. Martino},
{\it Blocks of restricted rational Cherednik algebras for $G(m,d,n)$},
J. of Algebra {\bf 397} (2014), 209-224.

\bibitem[Mos]{MoserOriginal} {\sc J. Moser},
{\it Three integrable Hamiltonian systems connected with isospectral deformations}, 
Advances in Math. {\bf 16} (1975), 197--220.

%

\bibitem[Nam1]{Namikawa-Poisson} {\sc Y. Namikawa},
{\it Poisson deformations of affine symplectic varieties},
Duke Math. J. {\bf 156}, 51--85 (2011).

\bibitem[Nam2]{Namikawa-Birational} {\sc Y. Namikawa},
{\it Poisson deformations and birational geometry},
J. Math. Sci. Univ. Tokyo {\bf 22} (2015), 339--359.

%
%
%
%
%

\bibitem[Sage]{Sage} {\sc The Sage Developers},
{S}ageMath, the {S}age {M}athematics {S}oftware {S}ystem, {\tt https://www.sagemath.org}.

\bibitem[Sch]{SchmittThesis} {\sc J. Schmitt},
{\it On $\mathbb{Q}$-factorial terminalizations of symplectic linear quotient singularities},
PhD thesis, University of Kaiserslautern--Landau (2023).

\bibitem[ShTo]{ST} {\sc G. C. Shephard \& J. A. Todd},
{\it Finite unitary reflection groups},
Canad. J. Math. {\bf 6} (1954), 274-304.


%

\bibitem[Thi3]{thiel-counter} {\sc U. Thiel},
{\it A counter-example to Martino’s conjecture about generic Calogero–Moser families},
Algebr. Represent. Theory {\bf 17} (5) (2014), 1323--1348.

\bibitem[Thi1]{champ} {\sc U. Thiel},
{\it CHAMP: A Cherednik Algebra Magma Package},
LMS J. Comput. Math. {\bf 18} (2015), 266--307.

\bibitem[Thi2]{thiel-blocks} {\sc U. Thiel},
{\it Blocks in flat families of finite-dimensional algebras},
Pac. J. Math. {\bf 295} (2018), 191--240.

%


\bibitem[Zas]{zaslavsky}{\sc T. Zaslavsky},
{\it Facing up to arrangements: face-count formulas for partitions of space by hyperplanes},
Mem. Amer. Math. Soc. {\bf 154} (1975).

\end{thebibliography}
\end{document}